\documentclass[11pt]{amsart}

\usepackage{latexsym}
\usepackage{amssymb}
\usepackage{amscd}
\usepackage[all]{xy}
\usepackage[mathscr]{euscript}
\usepackage{dsfont}

\normalsize

\RequirePackage{xspace}

\newcommand{\thought}[1]{}
\renewcommand{\thought}[1]{ \textbf{[#1]}}

\usepackage{enumerate}        
\newenvironment{roenumerate}{\begin{enumerate}[\upshape (i)]}{\end{enumerate}}

\newcommand\nc {\newcommand}
\newcommand\rnc{\renewcommand}

\newcount\blopone
\newcount\xone
\newcount\xtwo
\newcount\ytwo

\newtheorem{theorem}{Theorem}[section]
\newtheorem{prop}[theorem]{Proposition}
\newtheorem{com}[theorem]{Comment}
\newtheorem{apl}[theorem]{Application}
\newtheorem{exercise}[theorem]{Exercise}
\newtheorem{redu}[theorem]{Reduction}
\newtheorem{refinement}[theorem]{Refinement}
\newtheorem{summary}[theorem]{Summary}
\newtheorem{importnota}[theorem]{Important Notation}
\newtheorem{prblm}[theorem]{Problem}
\newtheorem{notation}[theorem]{Notation}
\newtheorem{explanation}[theorem]{Explanation}
\newtheorem{defin}[theorem]{Definition}
\newtheorem{caution}[theorem]{Caution}
\newtheorem{remark}[theorem]{Remark}
\newtheorem{reminder}[theorem]{Reminder}
\newtheorem{illustration}[theorem]{Illustration}
\newtheorem{observation}[theorem]{Observation}
\newtheorem{lemma}[theorem]{Lemma}
\newtheorem{construction}[theorem]{Construction}
\newtheorem{discussion}[theorem]{Discussion}
\newtheorem{corollary}[theorem]{Corollary}
\newtheorem{example}[theorem]{Example}
\newtheorem{conclusion}[theorem]{Conclusion}
\newtheorem{sketch}[theorem]{Sketch}
\newtheorem{triviality}[theorem]{Triviality}
\newtheorem{proto}[theorem]{Prototype Quasifibration}
\newtheorem{cauex}[theorem]{Cautionary Example}
\newtheorem{hypo}[theorem]{Hypothesis}
\newtheorem{subth}{ }[theorem]
\newtheorem{case}{Case}[theorem]
\newtheorem{ssubth}{ }[subth]
\newtheorem{facts}[theorem]{Facts}
\newtheorem{history}[theorem]{Historical Survey}
\newtheorem{proofs}[theorem]{Discussion of the Proofs, Old and New}

\nc\tri[1]{\begin{triviality}
\label{#1}}
\nc\fac[1]{\begin{facts}
\label{#1}
\begin{em}}
\nc\app[1]{\begin{apl}
\label{#1}
\begin{em}}
\nc\skt[1]{\begin{sketch}
\label{#1}
\begin{em}}
\nc\hst[1]{\begin{history}
\label{#1}
\begin{em}}
\nc\pfs[1]{\begin{proofs}
\label{#1}
\begin{em}}
\nc\cas[1]{\begin{case}
\label{#1}
\begin{em}}
\nc\rfn[1]{\begin{refinement}
\label{#1}}
\nc\prt[1]{\begin{proto}
\label{#1}}
\nc\lem[1]{\begin{lemma}
\label{#1}}
\nc\pro[1]{\begin{prop}
\label{#1}}
\nc\thm[1]{\begin{theorem}
\label{#1}}
\nc\dis[1]{\begin{discussion}
\label{#1}
\begin{em}}
\nc\cor[1]{\begin{corollary}
\label{#1}}
\nc\dfn[1]{\begin{defin}
\label{#1}}
\nc\sthm[1]{\begin{subth}
\label{#1}}
\nc\exm[1]{\begin{example}
\label{#1}
\begin{em}}
\nc\obs[1]{\begin{observation}
\label{#1}
\begin{em}}
\nc\plm[1]{\begin{prblm}
\label{#1}
\begin{em}}
\nc\rmk[1]{\begin{remark}
\label{#1}
\begin{em}}
\nc\rmd[1]{\begin{reminder}
\label{#1}
\begin{em}}
\nc\ntn[1]{\begin{notation}
\label{#1}
\begin{em}}
\nc\exe[1]{\begin{exercise}
\label{#1}
\begin{em}}
\nc\xpl[1]{\begin{explanation}
\label{#1}
\begin{em}}
\nc\smr[1]{\begin{summary}
\label{#1}
\begin{em}}
\nc\cau[1]{\begin{caution}
\label{#1}
\begin{em}}
\nc\hyp[1]{\begin{hypo}
\label{#1}
\begin{em}}
\nc\imn[1]{\begin{importnota}
\label{#1}
\begin{em}}
\nc\rdn[1]{\begin{redu}
\label{#1}
\begin{em}}
\nc\cax[1]{\begin{cauex}
\label{#1}
\begin{em}}
\nc\cmt[1]{\begin{com}
\label{#1}
\begin{em}}
\nc\con[1]{\begin{construction}
\label{#1}
\begin{em}}
\nc\ill[1]{\begin{illustration}
\label{#1}
\begin{em}}
\nc\ssthm[1]{\begin{ssubth}
\label{#1}
\begin{em}}
\nc\cnc[1]{\begin{conclusion}
\label{#1}
\begin{em}}

\nc\elem{\end{lemma}}
\nc\erdn{\end{em}\end{redu}}
\nc\erfn{\end{refinement}}
\nc\eprt{\end{proto}}
\nc\ethm{\end{theorem}}
\nc\ecor{\end{corollary}}
\nc\edfn{\end{defin}}
\nc\esthm{\end{subth}}
\nc\epro{\end{prop}}
\nc\etri{\end{triviality}}
\nc\eexm{\end{em}
\end{example}}
\nc\eobs{\end{em}
\end{observation}}
\nc\ecmt{\end{em}
\end{com}}
\nc\efac{\end{em}
\end{facts}}
\nc\eapp{\end{em}
\end{apl}}
\nc\ermk{\end{em}
\end{remark}}
\nc\ermd{\end{em}
\end{reminder}}
\nc\eill{\end{em}
\end{illustration}}
\nc\eplm{\end{em}
\end{prblm}}
\nc\ecas{\end{em}
\end{case}}
\nc\eskt{\end{em}
\end{sketch}}
\nc\ecau{\end{em}
\end{caution}}
\nc\ecax{\end{em}
\end{cauex}}
\nc\eimn{\end{em}
\end{importnota}}
\nc\entn{\end{em}
\end{notation}}
\nc\eexe{\end{em}
\end{exercise}}
\nc\expl{\end{em}
\end{explanation}}
\nc\edis{\end{em}
\end{discussion}}
\nc\econ{\end{em}
\end{construction}}
\nc\esmr{\end{em}
\end{summary}}
\nc\ehst{\end{em}
\end{history}}
\nc\epfs{\end{em}
\end{proofs}}
\nc\ehyp{
\end{em}
\end{hypo}}
\nc\ecnc{\end{em}
\end{conclusion}}
\nc\essthm{\end{em}
\end{ssubth}}

\nc\sst{\scriptstyle}
\newcommand{\comment}[1]{}
\newcommand{\ri}{\longrightarrow}
\newcommand{\sr}{\rightarrow}

\newcommand{\D}{{\mathbf D}}

\newcommand{\oo}{\otimes}

\nc\op{^{\hbox{\rm\tiny op}}}
\nc\mth{^{\hbox{\rm\tiny th}}}

\nc\script{\mathscr}
\nc\z{\zeta}
\nc\bc{{\mathbb{BC}}}
\nc\KK{{\mathbb{K}}}
\nc\ct{{\script T}}
\nc\cf{{\script F}}
\nc\cg{{\script G}}
\nc\ch{{\script H}}
\nc\ck{{\script K}}
\nc\cl{{\script L}}
\nc\cm{{\script M}}
\nc\cv{{\script V}}
\nc\ce{{\script E}}
\nc\cs{{\script S}}
\nc\car{{\script R}}
\nc\cd{{\script D}}
\nc\cc{{\script C}}
\nc\ca{{\script A}}
\nc\ci{{\script I}}
\nc\cj{{\script J}}
\nc\co{{\script O}}
\nc\cu{{\script U}}
\nc\cw{{\script W}}
\nc\cx{{\script X}}
\nc\Cp{{\script P}}
\nc\cq{{\script Q}}
\nc\cy{{\script Y}}
\nc\cz{{\script Z}}
\nc\bd{\begin{description}}
\nc\ed{\end{description}}
\nc\ctob{{\script C}at\big(\ci^{op},\ca\big)}
\nc\clim{{\ds\mathop{\rm lim}_{\ds\longleftarrow}}\,}
\nc\climi{\clim_{\!i}\,}
\nc\climn{\clim^{\!n}\,}
\nc\colim{{\ds\mathop{\rm colim}_{\ds\la}}}
\nc\colimj{{\ds\mathop{\rm colim}_{\ds\la}}{}_{j\,}}
\nc\oa{\overline{\ca}}
\nc\s{\sigma}
\nc\ta{\tau}
\nc\os{\overline\sigma}
\nc\ot{\overline\tau}
\nc\T{\Sigma}
\nc\Tm{\Sigma^{-1}}
\nc\de[1]{{\mathop{\rm deg(#1)}}}
\nc\Ad[1]{\mathop{\rm Ad}(#1)}
\nc\ad[1]{\mathop{\rm ad}(#1)}
\nc\kth{{\it K}--theory}
\nc\loc[1]{{\text{\rm Loc}(#1)}}
\nc\coloc[1]{{\text{\rm Coloc}(#1)}}

\def\der #1 {D\left(#1\right)}
\nc\prf{\begin{proof}}
\nc\eprf{\end{proof}}
\nc\ds{\displaystyle}
\nc\Tor{\text{\rm Tor}}

\nc\cb{{\script B}}
\nc\ab{{\script A}b}

\nc\be{\begin{roenumerate}}
\nc\ee{\end{roenumerate}}

\nc\cat[1]{{\script C}at\Big({\big\{#1\big\}}\op\,\,,\,\,\ab\Big)}
\nc\csab{{\script C}at\big(\cs^{op},\ab\big)}
\nc\ctab{{\script C}at\Big({\{\ct^\alpha\}}^{op},\ab\Big)}
\nc\csex{{\script E}x\big(\cs^{op},\ab\big)}
\nc\ctex{{\script E}x\Big({\{\ct^\alpha\}}^{op},\ab\Big)}
\nc\sub{\qquad\subset\qquad}
\nc\ctr[1]{{\left.\ct\left(-,#1\right)\right|}_{\cs}}
\nc\ctrf[2]{{\left.\ct\left(#1,#2\right)\right|}_{\cs}}
\nc\Ctr[1]{{\left.\ct\left(-,#1\right)\right|}_{\ct^\alpha}}
\nc\Ctrf[2]{{\left.\ct\left(#1,#2\right)\right|}_{\ct^\alpha}}

\nc\la{\longrightarrow}
\nc\nin{\noindent}
\nc\cad[1]{\text{card}(#1)}
\nc\eq{\quad=\quad}
\nc\BA{\begin{array}{c}}
\nc\EA{\end{array}}
\nc\barr{
\[
\begin{array}{cccccccccccccccc}
}
\nc\earr{
\end{array}
\]
}
\nc\as[1]{{\langle S\rangle}^{#1}}
\nc\sh{\text{\it shift}}

\nc\yy[1]{{\left.\ct\left(-,#1\right)\right|}_{\ct^c}}
\nc\vrep[2]{{\left.\ct\left(#1,#2\right)\right|}_{\ct^\alpha}}
\nc\da{\downarrow}
\nc\Hom{{\mathop{\rm Hom}}}
\nc\HHom{{\script H}{\mathop{\rm om}}}
\nc\End{{\mathop{\rm End}}}
\nc\Ext{{\mathop{\rm Ext}}}
\nc\mr\Modtc
\nc\PExt{{\mathop{\rm PExt}}}
\nc\stm{\text{\rm stmod}(kG)}
\nc\stM{\text{\rm StMod}(kG)}
\nc\e{\varepsilon}
\nc\p{\varphi}

\nc\rs{\s^{-1}A}
\nc\br{{\{\s^{-1}A\}}}
\nc\ra\ri
\nc\y[1]{\mathbf{y}#1}
\nc\x[1]{\mathbf{z}#1}
\nc\mmod[1]{\text{\rm mod--}#1}
\nc\Mod[1]{#1\text{--\rm Mod}}
\nc\Md {\ensuremath{\mathop{\textup{Mod}}}}
\rnc\mod[1]{\ensuremath{\mathop{#1\textup{--mod}}}\xspace}
\nc\MMod[1]{\text{Mod-}#1}
\nc\Modtc{\Mod{\ct^c}}
\nc\pgldim[1]{\mathop{\rm pgldim}\,#1}
\nc\tf{{\rm [TR5]}}
\nc\tfs{{\rm [TR5$^*$]}}
\nc\Fun{\text{\rm Funct}(F\op,\ab)}
\nc\sym{\text{\rm Sym}}
\nc\sgn{\text{\rm sgn}}
\nc\Pro{\text{\rm Prod}^{}_\alpha(F\op,\ab)}
\nc\Yt[1]{{\left.\Hom_\ct^{}\left(-,#1\right)\right|}_F^{}}
\nc\dl{\delta}
\nc\Proj[1]{#1\text{--\rm Proj}}
\nc\proj[1]{#1\text{--\rm proj}}
\nc\Flat[1]{#1\text{--\rm Flat}}
\nc\Inj[1]{#1\text{--\rm Inj}}
\nc\Ima{\mathrm{Im}}
\nc\Ker{\mathrm{Ker}}
\nc\ov{\overline}
\nc\wt{\widetilde}
\nc\wh{\widehat}
\nc\ph{\varphi}
\nc\tstr{{\it t}--structure}
\nc\tstrs{{\it t}--structures}
\nc\spec[1]{{\text{\rm Spec}\left(#1\right)}}

\nc\EProd{\text{\rm EProd}}
\nc\ECoprod{\text{\rm ECoprod}}
\nc\Prod{\text{\rm Prod}}
\nc\ldimp{\text{\rm LDim}^{\prod}}
\nc\ldimc{\text{\rm LDim}^{\coprod}}
\nc\gen[2]{{\langle#1\rangle}^{}_{#2}}
\nc\Gen[2]{{\big\langle#1\big\rangle}^{}_{#2}}
\nc\genu[3]{{\langle#1\rangle}^{[#3]}_{#2}}
\nc\ogen[1]{\ov{\langle#1\rangle}}
\nc\ogenun[2]{\ov{\langle#1\rangle}_{#2}^{}}
\nc\ogenu[3]{\ov{\langle#1\rangle}^{[#3]}_{#2}}
\nc\ogenul[3]{\ov{\langle#1\rangle}^{(-\infty,#3]}_{#2}}
\nc\ogenuf[3]{\ov{\langle#1\rangle}^{[#3,\infty)}_{#2}}
\nc\genuf[3]{{\langle#1\rangle}^{[#3,\infty)}_{#2}}
\nc\genul[3]{{\langle#1\rangle}^{(-\infty,#3]}_{#2}}
\nc\dperf[1]{\D^{\mathrm{perf}}(#1)}
\nc\dperfs[2]{\D_{#1}^{\mathrm{perf}}(#2)}
\nc\dcoh{\mathbf{D}^b_{\mathrm{coh}}}
\nc\dcohs[1]{\mathbf{D}^b_{\mathrm{coh},#1}}
\newcommand{\Dqc}{{\mathbf D_{\text{\bf qc}}}}
\newcommand{\Dqcs}[1]{{\mathbf D_{\text{\bf qc},#1}}}
\newcommand{\Dqcsmi}[1]{{\mathbf D^-_{\text{\bf qc},#1}}}
\newcommand{\Dqcspl}[1]{{\mathbf D^+_{\text{\bf qc},#1}}}
\newcommand{\Dqcsb}[1]{{\mathbf D^b_{\text{\bf qc},#1}}}

\nc\dmcoh{\mathbf{D}^-_{\mathrm{coh}}}
\nc\dmcohs[1]{\mathbf{D}^-_{\mathrm{coh},#1}}
\nc\dscoh{\mathbf{D}^{}_{\mathrm{coh}}}
\nc\RHHom{{\script{RH}}{\mathrm{om}}}
\nc\Coprod{\mathrm{Coprod}}
\nc\COprod{\mathrm{coprod}}
\nc\add{\mathrm{add}}
\nc\Add{\mathrm{Add}}
\nc\Smr{\mathrm{smd}}
\nc\id{\mathrm{id}}
\nc\LL{\mathbf{L}}
\nc\R{\mathbf{R}}
\nc\CC{\mathbf{C}}
\nc\wi{\wt{\text{\it\i}}}
\nc\exal{\ce\text{\it x}_\alpha(\ct^\alpha,\ab)}
\nc\exalz{\ce\text{\it x}_{\aleph_0}^{}(\ct^\alpha,\ab)}

\nc\tst[1]{\left({#1}^{\leq0},{#1}^{\geq0}\right)}
\nc\perf[1]{\script{P}\mathrm{erf}\left(#1\right)}
\nc\SF[1]{\mathrm{SF}\left(#1\right)}

\nc\one{\mathds{1}}
\nc\fc{\mathfrak{C}}
\nc\fl{\mathfrak{L}}
\nc\fs{\mathfrak{S}}
\nc\Prf{\text{\bf Perf}}
\nc\qc{\text{\bf qc}}
\nc\vect[1]{\script{V}\mathit{ect}\left(#1\right)}
\nc\vectgr[1]{\script{V}\mathit{ect}\script{G}\mathit{r}\left(#1\right)}
\nc\Ch{\mathbf{Ch}}

\nc\hoco{
\begin{picture}(40,10)
\put(20,0){\makebox(0,0)[b]{\text{\rm Hocolim}}}
\put(5,-2){\vector(1,0){30}}
\end{picture}\,\,}

\nc\holim{{
\begin{picture}(40,10)
\put(20,0){\makebox(0,0)[b]{\text{\rm Holim}}}
\put(35,-2){\vector(-1,0){30}}
\end{picture}}}

\renewcommand{\leq}{\leqslant}
\renewcommand{\geq}{\geqslant}

\begin{document}

\author{Amnon Neeman}\thanks{The research was supported 
  by the Australian Research Council
  (grants number DP200102537 and number DP210103397), and 
  by the Deutsche
  Forschungsgemeinschaft (SFB-TRR 358/1 2023 - 491392403).
  Also: part of the work was done at the 
  Isaac Newton Institute for Mathematical Sciences, Cambridge, with funding
  from  EPSRC grant number EP/R014604/1. This was part of
  the program \emph{K-theory, algebraic cycles and motivic homotopy theory.}
  The author
  would like to thank the Newton Institute for their
  support and hospitality
during this program. Another 
part of the work was done during the trimester program \emph{Spectral
Methods in Algebra, Geometry, and Topology} at the Hausdorff Institute
in Bonn. It is a pleasure to thank the Hausdorff Institute for its
hospitality, and acknowledge the funding
by the Deutsche Forschungsgemeinschaft (DFG) under Excellence
Strategy EXC-2047/1-390685813.}
\address{Dipartimento di Matematica ``F.\ Enriques''\\
        Universit{\`a} degli Studi di Milano\\
        Via Cesare Saldini 50\\
	20133 Milano\\
        ITALY}
\email{amnon.neeman@unimi.it}

\title{Bounded {\it t}--structures on the category of perfect complexes}

\begin{abstract}
Let $X$ be a finite-dimensional, noetherian scheme.
Antieau, Gepner and Heller conjectured that its
derived category of perfect complexes has a bounded \tstr\ 
if and only if $X$ is regular. 

We prove a generalization, and to do so we sharpen some
of the techniques so far obtained in
the theory of approximable triangulated categories.
\end{abstract}

\subjclass[2010]{Primary 14F08, secondary 18G80, 19D35}

\keywords{Derived categories, perfect complexes, 
bounded {\it t}--structures}

\maketitle

\tableofcontents

\setcounter{section}{-1}

\section{Introduction}
\label{S0}

The notion of \tstr\ originated in 
Be{\u\i}linson, Bernstein
and Deligne~\cite{BeiBerDel82}. Given a triangulated category $\ct$, a 
\tstr\ on $\ct$  is a pair of full subcategories 
$\big(\ct^{\leq0},\ct^{\geq0}\big)$ satisfying
axioms we will not repeat here; the reader can find them in 
Definition~\ref{D100.1}.
We recall that a \tstr\ $\big(\ct^{\leq0},\ct^{\geq0}\big)$
on the category $\ct$ is called \emph{bounded}
if, for every object $X\in\ct$, there exists an integer $n>0$
with $\T^nX\in\ct^{\leq0}$ and $\T^{-n}X\in\ct^{\geq0}$. And
a (bounded) \tstr\ on a stable infinity category $\cm$ is by
definition a (bounded) \tstr\ on its homotopy category.

Let $\cm$ be a stable infinity category. From
Antieau, Gepner and 
Heller~\cite[Theorems~1.1 and 1.2]{Antieau-Gepner-Heller19}
we learn that there are obstructions, in the negative \kth\ of
$\cm$, to the existence of bounded \tstr s on $\cm$.

This is an amazing result with many fascinating
implications---the reader is referred to
the survey~\cite{Neeman22B} for more.
For the current
article we confine ourselves to noting that singular, 
noetherian schemes $X$ often have nonvanishing
negative \kth\ and, as observed in 
\cite[Corollary~1.4]{Antieau-Gepner-Heller19},
this places restrictions on whether the category $\dperf X$
can have bounded \tstr s, and if so what kind of heart these \tstr s
are allowed to have.
Motivated by what they could prove,
Antieau, Gepner and 
Heller went on to formulate
\cite[Conjecture~1.5]{Antieau-Gepner-Heller19}: it
predicted that, if $X$ is a finite-dimensional, noetherian scheme, then 
the category $\dperf X$ has a bounded \tstr\ if and only if
$X$ is regular. In this article we prove a major generalization.

\thm{T0.1}
Let $X$ be a noetherian, finite-dimensional scheme, and let $Z\subset X$
be a closed subset. Let $\dperfs ZX$ be the category of perfect 
complexes on $X$ whose cohomology is supported on $Z$.

Then the category $\dperfs ZX$
has a bounded \tstr\ if and only if $Z$ is contained in the regular 
locus of $X$.
\ethm

Note that if $X$ is noetherian and finite-dimensional, and $Z$ is contained
in the regular locus of $X$, then 
$\dperfs ZX=\dcohs Z(X)$; that is every bounded complex of coherent sheaves,
with cohomology supported on $Z$, is a perfect complex.
And $\dcohs Z(X)$ always
has an obvious bounded \tstr---the standard \tstr\ on $\Dqcs Z(X)$ restricts
to a bounded
\tstr\ on the full subcategory $\dcohs Z(X)$. What needs proof is that
the existence of a bounded \tstr\ on
$\dperfs ZX$ implies regularity.

\rmk{R0.3}
In the light of the evidence available at the time, Antieau, Gepner 
and Heller~\cite[Conjecture~1.5]{Antieau-Gepner-Heller19} was bold
and daring. 
After all: there are many singular schemes with vanishing
negative \kth, for example all singular zero-dimensional schemes. 
For such singular schemes,
\cite[Theorems~1.1 and 1.2]{Antieau-Gepner-Heller19}
place no restrictions on any potential bounded \tstr.
Thus 
the non-existence of bounded \tstr s, for all singular schemes, goes 
way beyond what the $K$-theoretic obstructions proved.

The best conjectures are the courageous ones that turn out to be true.

And Theorem~\ref{T0.1} is much more general yet. Thus negative \kth\ most 
definitely isn't the only obstruction to the existence of bounded \tstr s,
and it becomes interesting to figure out what other obstructions
there are, in a generality that goes beyond $\dperfs ZX$. This 
should be studied.

The current article does not proceed that way;
as we will sketch in the coming paragraphs, the proof we give of 
Theorem~\ref{T0.1} is not obstruction-theoretic. 
\ermk

The key to the proof is a certain equivalence relation on \tstr s, introduced
by the author in~\cite{Neeman17A}. 
We will give a self-contained treatment, reminding 
the reader of this equivalence 
relation and its
origins in Section~\ref{S9},
more precisely in the discussion leading up to
Definition~\ref{R9.13}. The strategy of the proof
of Theorem~\ref{T0.1}
will follow two steps.
\be
\item
Assume $\big({\dperfs ZX}^{\leqslant0},{\dperfs ZX}^{\geqslant0}\big)$ 
is a bounded \tstr\ on
$\dperfs ZX$. Then the \tstr\ on $\Dqcs Z(X)$ generated by
${\dperfs ZX}^{\leqslant0}$, in the sense of 
Alonso, Jerem{\'{\i}}as and Souto~\cite{Alonso-Jeremias-Souto03},
is equivalent to the standard \tstr\ on $\Dqcs Z(X)$.
\item
We then use a theorem proving that all objects
of $\dcohs Z(X)$ have arbitrarily good approximations
by objects of $\dperfs ZX$, see Theorem~\ref{T29.1}. 
\ee

If $Z=X$ then it is customary to drop the $Z$ in
the notation, as in $\Dqcs Z(X)=\Dqc(X)$, $\dperfs ZX=\dperf X$ and 
$\dcohs Z(X)=\dcoh(X)$. And for $Z=X$ the 
result (ii) is relatively old, it can be found 
already in Lipman and 
Neeman~\cite[Theorem~4.1]{Lipman-Neeman07}. Assuming furthermore
that
$Z=X$ is separated, there is a trick, explained in
Section~\ref{S27}, allowing us to deduce (i) 
from the theory developed in~\cite{Neeman17A}.
This was the approach taken in an earlier
incarnation of the present article: the resulting paper
was much shorter, but the theorems weren't nearly as
sharp. 

Part of the interest, of the 
current, improved version of the article, is 
the generalization of the requisite
machinery to cover the case of $\Dqcs Z(X)$; see
Theorem~\ref{T27.1} for (i) and Theorem~\ref{T29.1} for (ii). 
Fortunately for us, the work 
for proving  Theorem~\ref{T29.1} has already been done for
us, in the Stacks Project. See the references
 in Section~\ref{S29}, in the
paragraph just preceding
the precise statement of Theorem~\ref{T29.1}. 
Thus the real technical improvement that is new to this
article is Theorem~\ref{T27.1}. The short summary, in
the language of~\cite{Neeman17A}, is that we
prove (a) that the category $\Dqcs Z(X)$ is weakly approximable,
and (b) that the standard \tstr\ on $\Dqcs Z(X)$ is in
the preferred equivalence class.

Because the equivalence relation on \tstr s is essential to
the article, we give a self-contained, complete treatment of
the proof of (b). The assertion (a) is proved in the (separate)
Section~\ref{S28}, which is designed so that the reader can
safely skip it. The reader who does choose to read Section~\ref{S28}
will be assumed to have some familiarity with the
author's~\cite{Neeman17A}, we allow ourselves to
freely cite results
from there.
And in Section~\ref{S400}
there will be further discussion, of the combination of 
Theorem~\ref{T27.1} with Theorem~\ref{T29.1} and how it fits
in with the theory of approximability in triangulated categories.

It follows from (i), coupled with
the easy Lemma~\ref{L3.1}, 
that all bounded \tstr s on $\dperfs ZX$ are equivalent.
If $Z$ meets the singular locus of $X$, then Theorem~\ref{T0.1} says that
there are no bounded \tstr s on $\dperfs ZX$; this
makes the assertion that the bounded \tstr s are 
all equivalent true but empty.
It is therefore natural
to ask what is true about $\dcohs Z(X)$---which has at 
least one bounded \tstr\ for every $X$. The last pair of results
in the article, Theorems~\ref{T31.3} and \ref{T3.3}, give
conditions under
which we can prove that
the bounded \tstr s on $\dcohs Z(X)$ are all equivalent.
A summary of what we know is presented in Remark~\ref{R997.3}.

It seems eminently plausible that the results about
$\dcohs Z(X)$ are not best
possible. They are all we can prove at the moment.

We should mention that, until now, there were few known results in
the direction of 
Antieau, Gepner and Heller \cite[Conjecture~1.5]{Antieau-Gepner-Heller19}.
There were of course the {\it K}-theoretic theorems
of~\cite{Antieau-Gepner-Heller19} which motivated the 
conjecture---the reader can find an excellent account of
this background in the introduction 
to~\cite{Antieau-Gepner-Heller19}.
Harry Smith~\cite{SmithHarry19} proved the conjecture for affine $X$, and 
Gepner and Haesemeyer studied the class of bounded \tstr s
that behave well with respect to restriction to affines.
For those \tstr s, reduction to Harry Smith's result can be 
applied.

\medskip

\nin
{\bf Acknowledgements.}\ \ The author would like to thank
Christian Haesemeyer, who kindly presented an outline of his joint  project 
with Gepner. This is what 
spurred the author to try the approach given here. The author is 
also grateful
to Sasha Kuznetsov and Evgeny Shinder,
whose tremendously helpful and insightful
comments much clarified earlier versions. And I'm grateful to 
Jack Hall for many corrections, and for pointing me 
to the proof of Theorem~\ref{T29.1} in
the Stacks Project.
Further thanks go
to Paul Balmer, for corrections and 
improvements on earlier drafts, and to Alberto Canonaco and Paolo
Stellari, whose questions led me to improve an earlier version
of the article by going to the effort of proving the weak
approximability to $\Dqcs Z(X)$. Let me not explain this in
detail here,
but it turns out that this refinement, of an earlier incarnation
of Theorem~\ref{T27.1}, is useful in addressing a question
that Canonaco and Stellari asked. This will appear as a
separate article.

And finally let me thank two anonymous referees for their careful,
thorough reading of the manuscript, and for many
valuable suggestions of
improvements to the organization of the paper.

\section{Notational conventions}
\label{S100}

Let $\ct$ be a triangulated category. Following topologists, the \emph{shift} or
\emph{suspension} functor
will be denoted $\T:\ct\la\ct$. Thus we write $\T X$ for what's often
written 
$X[1]$.

Let $\ct$ be a
triangulated category and let $\ca\subset\ct$ be any collection of objects.
Then
$\ca^\perp$ denotes the full subcategory of all objects annihilated by
$\Hom(\ca,-)$, while $^\perp\ca$ denotes the
full subcategory of all objects annihilated by
$\Hom(-,\ca)$. And the reader should note that, since $\T$ is an automorphism
of the category, we have $(\T\ca)^\perp=\T(\ca^\perp)$ and
$^\perp(\T\ca)=\T(^\perp\ca)$. We will therefore
usually permit ourselves to omit
the brackets when performing these operations.

Our \tstr s will be denoted cohomologically. Thus a \tstr\ on a triangulated
category $\ct$ is a pair of full subcategories 
$(\ct^{\leqslant0},\ct^{\geqslant0})$ satisfying
the axioms of Be{\u\i}linson, Bernstein and 
Deligne~\cite[D\'efinition~1.3.1]{BeiBerDel82}. For the reader's convenience
we recall this:

\dfn{D100.1}
A \tstr\ on a triangulated category is a pair of full subcategories
$\tst\ct$ satisfying:
\be
\item
$\T\ct^{\leq0}\subset\ct^{\leq0}$ and $\ct^{\geq0}\subset\T\ct^{\geq0}$.
\item
$\T\ct^{\leq0}\subset{^\perp\ct^{\geq0}}$.
\item
For every object $B\in\ct$ there exists a triangle $A\la B\la C\la\T A$ with
$A\in\T\ct^{\leq0}$ and with $C\in\ct^{\geq0}$.  
\ee
\edfn

\rmd{R100.3}
The expert will note that there is some redundancy in this old definition,
there exist more succinct formulations. For us it will be important to
remember that the containment in Definition~\ref{D100.1}(ii) is actually
an equality: for any \tstr\ one can prove
\[
\T\ct^{\leq0}={^\perp\ct^{\geq0}}\qquad\text{ and }\qquad\left(\T\ct^{\leq0}\right)^\perp=\ct^{\geq0}
\]
For any collection $\ca\subset\ct$ we have that $^\perp\ca$ is closed in
$\ct$ under extensions, direct summands and coproducts
(meaning those coproducts that
exist in $\ct$).
Dually $\ca^\perp$ is closed in $\ct$ under extensions, direct summands
and products. From the equalities above we deduce
\be
\item
Let $\ct$ be a triangulated category and let $\tst\ct$ be a \tstr\ on
it. Then, in addition to the inclusion $\T\ct^{\leq0}\subset\ct^{\leq0}$
of Definition~\ref{D100.1}(i), we have that
the subcategory $\ct^{\leq0}={^\perp\Tm\ct^{\geq0}}$
is closed in $\ct$ under extensions,  
direct summands and those coproducts that exist in $\ct$. Following
Keller and Vossieck~\cite{Keller-Vossieck88}, any
subcategory $\ca\subset\ct$ satisfying these closure properties
is called a \emph{pre-aisle.} If there is a \tstr\ with
$\ca=\ct^{\leq0}$
then the pre-aisle $\ca$ is called an \emph{aisle.}
\item
Let $\ct$ be a triangulated category and let $\tst\ct$ be a \tstr\ on
it. Then, in addition to the inclusion $\Tm\ct^{\geq0}\subset\ct^{\geq0}$
of Definition~\ref{D100.1}(i), we have that the
subcategory $\ct^{\geq0}=(\T\ct^{\leq0})^\perp$
is closed in $\ct$ under extensions,  
direct summands and those products that exist in $\ct$. Any
subcategory $\cb\subset\ct$ satisfying these closure properties
is called a \emph{pre-coaisle.} If there is a \tstr\
with $\cb=\ct^{\geq0}$
then the pre-coaisle $\cb$ is called a \emph{coaisle.}
\ee
\ermd

Following the standard (cohomological) notation, we set
\[
\ct^{\leq n}=\T^{-n}\ct^{\leq0}\qquad\text{ and }\qquad\ct^{\geq n}=\T^{-n}\ct^{\geq0}\ .
\]
And we remind the reader that the triangle, whose
existence is guaranteed in Definition~\ref{D100.1}(iii), is
unique up to canonical isomorphism and functorial. It is usually
written as
\[\xymatrix{
B^{\leq-1}\ar[r] & B\ar[r] & B^{\geq0}\ar[r] &
\T(B^{\leq-1})\ .
}\]

\section{A reminder of a few results from
the theory of approximable
triangulated categories}
\label{S9}

In this article we will use a few results from
the series of articles 
\cite{Neeman17,Neeman17A,Burke-Neeman-Pauwels18,Neeman18,Neeman18A}.
To help the reader unfamiliar with this
work, this section is a review of the small portion
of the general theory necessary for the current article.

\rmd{R9.1}
Let $A$ be a set of objects in a triangulated category $\ct$. We 
begin by recalling
the full subcategories $\COprod(A)$ and
$\Coprod(A)$ of $\ct$, whose definitions and elementary
properties may be found in \cite[Section~1]{Neeman17}. 

By combining
\cite[Definition~1.3(iii) and Observation~1.5]{Neeman17} we have
\be
\item
$\COprod(A)\subset\ct$ can be defined to be the smallest full
subcategory $\cs\subset\ct$
containing $A$ and closed under finite coproducts
and extensions. For the reader preferring symbols
and happy with the notation of \cite[Reminder~1.1]{Neeman17}, this
can be reformulated as saying that $\COprod(A)$ is the minimal 
$\cs\subset\ct$ satisfying
\[
A\subset\cs,\qquad\add(\cs)\subset\cs\qquad\text{and}\qquad\cs\star\cs\subset\cs.
\]  
\setcounter{enumiv}{\value{enumi}}
\ee
Now assume $\ct$ has coproducts. Then 
\cite[Definition~1.3(iv)]{Neeman17} introduces the ``big'' version, that is
it defines
\be
\setcounter{enumi}{\value{enumiv}}
\item
$\Coprod(A)\subset\ct$ is the smallest full
subcategory $\cs\subset\ct$
containing $A$ and closed under all coproducts
and extensions. In symbols
and using the notation of \cite[Reminder~1.1]{Neeman17}, we can 
restate this by saying that $\Coprod(A)$ is the minimal $\cs\subset\ct$ satisfying
\[
A\subset\cs,\qquad\Add(\cs)\subset\cs\qquad\text{and}\qquad\cs\star\cs\subset\cs.
\] 
\setcounter{enumiv}{\value{enumi}} 
\ee
\ermd

Next we recall the little

\lem{L9.3}
Assume that $\ct$ is a triangulated category with
coproducts, and let $A\subset\ct$ be a set of objects
 satisfying $\T A\subset A$.
Then, in the terminology of 
Alonso, Jerem{\'{\i}}as and 
Souto~\cite[Section~1]{Alonso-Jeremias-Souto03},
the subcategory $\Coprod(A)\subset\ct$ is the cocomplete
pre-aisle generated by $A$.
\elem

\prf
By definition we have that $A$ is contained in
$\Coprod(A)$, that $\Add\big(\Coprod(A)\big)\subset\Coprod(A)$,
that $\Coprod(A)\star\Coprod(A)\subset\Coprod(A)$, and that
$\Coprod(A)$ is the minimal
subcategory with these properties. It remains to
prove that $\Coprod(A)$ is closed under suspension. Here
our assumption that $\T A\subset A$ comes in; it gives the inclusion in
\[
\T\Coprod(A)\eq\Coprod(\T A)\sub\Coprod(A)\ .
\]
The equality is because the functor $\T$ is an automorphism respecting
extensions and (of course) coproducts.
And the conclusion is that $\Coprod(A)$ is closed under suspension.
\eprf

\rmk{R9.4}
We remind the reader that, as long as $\T A\subset A$ and therefore
$\T\Coprod(A)\subset\Coprod(A)$, it's automatic
that $\Coprod(A)$ is closed under direct summands. See for
example the Eilenberg swindle argument in
\cite[footnote on page~227]{Alonso-Jeremias-Souto00}.
\ermk

\rmd{R9.5}
As observed in 
\cite[Problem following Corollary~1.4]{Alonso-Jeremias-Souto03},
it becomes interesting to find sufficient conditions,
guaranteeing that the cocomplete pre-aisle $\Coprod(A)$ is 
an aisle---meaning conditions that guarantee the
existence of a (unique) \tstr\ $(\ct^{\leqslant0},\ct^{\geqslant0})$,
on the category $\ct$, with 
\[
\ct^{\leqslant0}=\Coprod(A).
\] 
There is an extensive literature on this,
but for our limited purposes the old result
in \cite[Theorem~A.1]{Alonso-Jeremias-Souto03} 
suffices. It asserts:
\be
\item
Let $\ct$ be a triangulated category with coproducts. Let
$A$ be a set of \emph{compact} objects in $\ct$ satisfying
$\T A\subset A$. Then $\Coprod(A)$ is an aisle.
\ee
The \tstr\ $(\ct^{\leqslant0},\ct^{\geqslant0})$
with $\ct^{\leqslant0}=\Coprod(A)$ is called the \tstr\ \emph{generated by $A$},
and the class of such \tstr s is referred to as the
\emph{compactly generated \tstr s.}
\ermd 

Later in the article we will have occasion to use the following
little result.

\lem{L3.1}
Let $\ct$ be a triangulated category with coproducts, and let $\cs\subset\ct$
be a full, triangulated subcategory. Assume we are given a \tstr\ 
$(\cs^{\leqslant0},\cs^{\geqslant0})$ on $\cs$, and suppose that $\cs^{\leqslant0}$ generates
a \tstr\ on $\ct$, in the sense of Reminder~\ref{R9.5}.
To spell it out: this means that 
$\big(\Coprod(\cs^{\leqslant0}),\T\Coprod(\cs^{\leqslant0})^\perp\big)$ is a 
\tstr\ on $\ct$.

Then the restriction to $\cs$ of the \tstr\ generated in $\ct$ by
$\cs^{\leqslant0}$ is equal to $(\cs^{\leqslant0},\cs^{\geqslant0})$.
This means that
\[\cs^{\leqslant0}=\cs\cap\Coprod(\cs^{\leqslant0})\qquad\text{ and }\qquad
\cs^{\geqslant0}=\cs\cap\T\Coprod(\cs^{\leqslant0})^\perp.
\]
and that, if $B\in\cs$ is an object, then the triangle
$B^{\leq-1}\la B\la B^{\geq0}\la\T(B^{\leq-1})$ is identical, whether we
form it in $\cs$ or in $\ct$.
\elem

\prf
If $B\in\cs$ is any object, then we may truncate it with
respect to the \tstr\ $\tst\cs$ to obtain
a distinguished triangle $B^{\leq-1}\la B\la B^{\geq0}\la\T(B^{\leq-1})$
with
\[
\begin{array}{ccccl}
B^{\leq-1}&\quad\in\quad& \cs^{\leqslant-1} &\sub& \Coprod(\cs^{\leqslant-1}) \\
B^{\geq0}&\quad\in\quad& \cs^{\geqslant0} &\sub& \Coprod(\cs^{\leqslant-1})^\perp\ .
\end{array}
\]
This means that the triangle
$B^{\leq-1}\la B\la B^{\geq0}\la\T(B^{\leq-1})$ is also the unique triangle
in $\ct$, giving the truncation with respect to the \tstr\ 
$\big(\Coprod(\cs^{\leqslant0}),\T\Coprod(\cs^{\leqslant0})^\perp\big)$.
The equalities
\[\cs^{\leqslant0}=\cs\cap\Coprod(\cs^{\leqslant0})\qquad\text{ and }\qquad
\cs^{\geqslant0}=\cs\cap\T\Coprod(\cs^{\leqslant0})^\perp.
\]
follow immediately.
\eprf

\rmd{R9.7}
One special case of Reminder~\ref{R9.5},
that will interest us in this article, is where
$G$ is a compact object of $\ct$, and $A=G(-\infty,0]$ is
given by
\[
G(-\infty,0]\eq\left\{
\T^i G\mid i\geqslant0
\right\}\ .
\]
In an abuse of the language we have used so far, we will call the 
\tstr\ generated by $G(-\infty,0]$ the \tstr\ generated by $G$,
and we will adopt the shorthand notation
\[
\ct_G^{\leqslant0}=\Coprod\big(G(-\infty,0]\big),\qquad
\ct_G^{\geqslant0}=\T\Coprod\big(G(-\infty,0]\big)^\perp.
\]
More generally: for any object $H\in\ct$ and any pair of integers 
$m\leqslant n$, possibly infinite,
\cite[Notation~1.11]{Neeman17} adopts the convention that
\[
H[m,n]\eq\left\{
\T^i H\mid -n\leqslant i\leqslant -m
\right\}\ .
\]
This set of objects is guaranteed to be stable under 
suspension only when 
$m=-\infty$.
\ermd

This leads us to the little

\lem{L9.9}
Let $G\in\ct$ be a compact generator, and let $H\in\ct$ be any compact
object. Then there exists an integer $n>0$ such that 
$H\in\ct_G^{\leqslant n}$.
\elem

\prf
In this proof we adopt the notation of
Bondal and Van den Bergh~\cite[beginning of 2.2]{BondalvandenBergh04},
where for any full subcategory $\ca\subset\ct$ we let $\Smr(\ca)$
be the full
subcategory of all direct summands of objects in $\ca$.

Because $G$ is a compact generator we have
$\ct=\Coprod\big(G(-\infty,\infty)\big)$, with the notation as in
Reminder~\ref{R9.7}. Thus $H\in\Coprod\big(G(-\infty,\infty)\big)$, but
by assumption $H\in\ct^c$. 
By \cite[Proposition~1.9(ii)]{Neeman17} we deduce that 
$H\in\ct^c\cap\Coprod\big(G(-\infty,\infty)\big)=
\Smr\Big(\COprod\big(G(-\infty,\infty)\big)\Big)$. 
But as $G(-\infty,\infty)=\bigcup_{n=1}^\infty G(-\infty,n]$, we have that
\cite[Definition~1.3(iii) and Observation~1.5]{Neeman17}
combine to yield
\[
\Smr\Big(\COprod\big(G(-\infty,\infty)\big)\Big)\eq\bigcup_{n=1}^\infty
\Smr\Big(\COprod\big(G(-\infty,n]\big)\Big)\ .
\]
Therefore there exists an integer $n>0$ with
\[
H\quad\in\quad
\Smr\Big(\COprod\big(G(-\infty,n]\big)\Big)
\sub
\Coprod\big(G(-\infty,n]\big)\eq
\ct_G^{\leqslant n}\ .
\]
\eprf

Let $G$, $H$ and $n$ be as in Lemma~\ref{L9.9}. Lemma~\ref{L9.9} tells
us that $H\in\ct_G^{\leqslant n}$,
hence $\T^mH\in\ct_G^{\leqslant n}$ for all $m\geqslant0$, that
is $H(-\infty,0]\subset\ct_G^{\leqslant n}$. Therefore
$\ct_H^{\leqslant0}\subset\ct_G^{\leqslant n}$,
and by symmetry we deduce

\cor{C9.11}
If $G$ and $H$ are both compact generators for $\ct$, then there
exists an integer $n>0$ such that
\[
\ct_G^{\leqslant-n}\sub\ct_H^{\leqslant0}\sub\ct_G^{\leqslant n}\ .
\]
\ecor

The following now becomes natural.

\dfn{R9.13}
Let $\ct$ be a triangulated category and let $(\ct_1^{\leqslant0},\ct_1^{\geqslant0})$ and
$(\ct_2^{\leqslant0},\ct_2^{\geqslant0})$ be two \tstr s on $\ct$. We declare these
\tstr s to be \emph{equivalent} if there exists an integer $n>0$ with
\[
\ct_1^{\leqslant-n}\subset\ct_2^{\leqslant0}\subset\ct_1^{\leqslant n}\ .
\]
\edfn

\rmk{R9.14}
Taking orthogonals, it is clear that the
$(\ct_1^{\leqslant0},\ct_1^{\geqslant0})$ and
$(\ct_2^{\leqslant0},\ct_2^{\geqslant0})$ are equivalent if and only if
there exists an integer $n$ with
\[
\ct_1^{\geqslant  n}\subset\ct_2^{\geqslant0}\subset\ct_1^{\geqslant-n}\ .
\]
\ermk

\rmd{R9.15}
It's obvious that equivalence of \tstr s is an equivalence relation. 
Corollary~\ref{C9.11} asserts that, if $G$ and $H$ are compact generators,
then the \tstr s generated by $G$ and by $H$ are equivalent. Thus
if the category $\ct$ has a compact generator then there is a
\emph{preferred equivalence class} of \tstr s, defined to be all those
equivalent to the one generated by $G$ for some compact generator $G$,
and hence for every compact generator $G$.
\ermd

\section{The proof of Theorem~\protect{\ref{T0.1}}}
\label{S29}

We start this section with a couple
of old theorems, but where the current generality
is recent.
To state the theorems we first recall some classical notation.

\rmd{R9.17}
Let $X$ be a
scheme and $Z\subset X$ a closed subset. The category $\Dqc(X)$ 
has for its  objects
the cochain complexes of $\co_X^{}$-modules with
quasicoherent cohomology. The full subcategory
$\Dqcs Z(X)\subset\Dqc(X)$ has for objects the complexes 
further restricted to be
acyclic on the open set $X-Z$.
\ermd

In sections~\ref{S27} and~\ref{S28} we will give
a self-contained proof of the following
theorem, as well as a discussion of its history.

\thm{T27.1}
Let $X$ be a quasicompact, quasiseparated scheme. Assume $Z\subset X$
is a closed subset such that $X-Z$ is quasicompact. 
Let $\Dqcs{Z}(X)\subset\Dqc(X)$ be as in Reminder~\ref{R9.17}.

Then the following assertions are true:
\be
\item
The compact objects in $\Dqcs{Z}(X)$ are those perfect complexes on
$X$ whose cohomology is supported in $Z$. We will denote the 
full subcategory of such objects $\dperfs{Z}X$.
\item
There exists a single compact object $G\in\Dqcs{Z}(X)$ generating 
$\Dqcs{Z}(X)$.
\item
The standard \tstr\ on $\Dqcs{Z}(X)$ is in the preferred equivalence
class. 
\item
The category $\Dqcs Z(X)$ is weakly approximable, as in
\cite[Definition~0.21]{Neeman17A}.
\ee
\ethm

The next theorem also has a long history,
but the sharp, recent version we need  is fortunately
already in the literature.
To state this
result, in full, glorious generality,  we remind
the reader of more classical terminology.

Let $X$ be a quasicompact, quasiseparated scheme. A complex $A\in\Dqc(X)$
is \emph{pseudocoherent}
if, for any open immersion 
$j:\spec R\la X$, we have that
the object $j^*A\in\Dqc\big(\spec R\big)\cong\D(R)$ is
isomorphic to an object in $\D^-(\proj R)$. To say it in words rather than 
symbols: for the object $A$ to qualify as pseudocoherent, the restriction
cochain
complexes $j^*A$, for every open immersion
$j:\spec R\la X$, must all
have bounded-above resolutions by finite-rank 
vector bundles.

The reader is referred to \cite[Definition~ 36.14.1 and 
Theorem~36.14.6]{stacks-project} for the statement and proof of
a result that is technically superior to Theorem~\ref{T29.1} below,
in a way
irrelevant to us here.

\thm{T29.1}
Let $X$ be a quasicompact, quasiseparated scheme, and assume $Z\subset X$
is a closed subset with quasicompact complement. 

Let $n>0$ be any positive integer.
Then every pseudocoherent object $F\in\Dqcs{Z}(X)$ admits a triangle
$E\la F\la D\la\T E$, with $E$ a compact object in $\Dqcs{Z}(X)$ and 
with $D\in\Dqcs{Z}(X)^{\leq-n}$.
\ethm

This section sets out to prove Theorem~\ref{T0.1}, and in
Theorem~\ref{T0.1} the scheme $X$ is assumed noetherian.
Hence all open subsets of $X$ are quasicompact and,
since asserting that $X$ is
quasiseparated is equivalent to the quasicompactness of
certain of its open subsets, it's also automatic
that $X$ is quasiseparated. We can
therefore freely apply Theorem~\ref{T29.1} and
all parts of
Theorem~\ref{T27.1} to $\Dqcs Z(X)$, for any closed subset $Z\subset X$.

We remark here that everything except Theorem~\ref{T27.1}(iv)
will in fact be used. The proof of Theorem~\ref{T27.1}(iv) is
therefore relegated to its own section~\ref{S28}, allowing the
reader to easily skip it.

Theorem~\ref{T0.1} is about bounded \tstr s on
$\dperfs ZX$, and the idea of
the proof will be to relate these to the \tstr s they generate
in $\Dqcs Z(X)$, in the sense of Reminder~\ref{R9.5}. A
key step in the proof is 
the following lemma. 

\lem{L1.3}
Suppose $X$ is a finite-dimensional, noetherian scheme and 
$Z\subset X$ is a closed 
subset.
Assume that $\big({\dperfs ZX}^{\leqslant0},{\dperfs ZX}^{\geqslant0}\big)$ is a bounded 
\tstr\ on $\dperfs ZX$.

Then the \tstr\ on $\Dqcs Z(X)$ generated by 
${\dperfs ZX}^{\leqslant0}$
is in the preferred equivalence class.
\elem

We will prove Lemma~\ref{L1.3} at the end of this section.
First we explain
how to deduce Theorem~\ref{T0.1} from the lemma.

\medskip

\nin
\emph{Proof that Theorem~\ref{T0.1} follows from Lemma~\ref{L1.3}.}\ \ 
The hypotheses on $Z\subset X$, made in Theorem~\ref{T0.1}, are as in
Lemma~\ref{L1.3}. Let
$\big({\dperfs ZX}^{\leqslant0},{\dperfs ZX}^{\geqslant0}\big)$ be a bounded 
\tstr\ on $\dperfs ZX$.
By the conclusion of Lemma~\ref{L1.3}
the \tstr\ on
$\Dqcs Z(X)$, generated by ${\dperfs ZX}^{\leqslant0}$,
is in the preferred equivalence class.
By Theorem~\ref{T27.1}(iii) it is equivalent 
to the standard \tstr\ on $\Dqcs Z(X)$.

And now it's time to come to the point: we will prove the regularity
in $X$ of every point in $Z\subset X$ by showing that the inclusion 
$\dperfs ZX\subset\dcohs Z(X)$ is an
equality. Choose therefore any object $F\in\dcohs Z(X)$, and we will
prove that $F$ must be a perfect complex. Shifting $F$ if necessary, we
may assume that $F\in\Dqcs Z(X)^{\geqslant0}$.

Put $\ca={\dperfs ZX}^{\leqslant0}$. The 
equivalence of the standard \tstr\ on $\Dqcs Z(X)$ with the \tstr\ 
generated by $\ca$ allows us to choose an
integer $n>0$ with
\[
\Dqcs Z(X)^{\leqslant -n}\sub\Coprod(\ca)\sub\Dqcs Z(X)^{\leqslant n}\ .
\]
Next we appeal to Theorem~\ref{T29.1}. It tells us that,
for any integer $M>0$, we can find an exact triangle $D\la E\la F$ with 
$E\in\dperfs ZX$ and with $D\in\Dqcs Z(X)^{\leqslant-M}$. Now 
put $M=1+2n$, with $n$ as above. In the triangle
$D\la E\la F$ of Theorem~\ref{T29.1} we have that 
\[
\begin{array}{ccccc}
D&\quad\in\quad&\Dqcs Z(X)^{\leqslant -1-2n}&\sub&\T^{1+n}\Coprod(\ca)\\
F&\quad\in\quad&\Dqcs Z(X)^{\geqslant0}&\sub&\T^{1+n}\Coprod(\ca)^\perp
\end{array}
\]
It means that, with respect to the \tstr\
$\big(\Coprod(\ca),\T\Coprod(\ca)^\perp\big)$ on $\Dqcs Z(X)$,
we have $D=E^{\leq-n-1}$ and $F=E^{\geq-n}$.
But now we appeal to Lemma~\ref{L3.1}; since the object
$E$ belongs to $\dperfs ZX$ its truncations in the larger
category $\Dqcs Z(X)$, with respect to the \tstr\ generated
by $\ca$, must belong to $\dperfs ZX$.
In particular $F\cong E^{\geq-n}\in\dperfs ZX$.
\hfill{$\Box$}

It remains to prove Lemma~\ref{L1.3}, but
before launching into the proof of  we give one more little
lemma.

\lem{L2.1}
Let $X$ be a noetherian scheme of finite dimension 
$m\geqslant0$. If $E\in\dperf X$
is a perfect complex on $X$, where $E\in\Dqc(X)^{\geqslant0}$ with respect
to the standard \tstr, then $E^\vee=\RHHom(E,\co_X^{})$ belongs
to $\Dqc(X)^{\leqslant m}$.
\elem

\prf
The assertion is local in $X$, hence we may assume $X=\spec R$ with $R$ 
a noetherian, local ring of Krull dimension $\leqslant m$. Under the equivalence
of categories $\dperf X\cong \D^b(\proj R)$, we have that $E$ 
corresponds to a bounded complex of finitely generated free
modules 
\[\xymatrix{
0\ar[r] & P^{-k}\ar[r] & P^{-k+1}\ar[r] &\cdots\ar[r] &
P^{\ell-1}\ar[r] & P^{\ell}\ar[r] & 0\ .
}\] 
And since $E$ is assumed to belong to $\Dqc(X)^{\geqslant0}$, it follows that
this complex has its cohomology in degrees $\geqslant0$. 

Therefore the truncated complex
\[\xymatrix{
0\ar[r] & P^{-k}\ar[r] & P^{-k+1}\ar[r] &\cdots\ar[r] &
P^{-1}\ar[r] & P^{0}\ar[r] & 0
}\]
has its cohomology concentrated in degree 0; this complex is a free
resolution of a finite $R$--module $M$. But now the 
Auslander-Buchsbaum formula tells us that 
\[
\text{proj.~dim}(M)+\text{depth}(M)\eq\text{depth}(R)\quad\leqslant\quad m,
\]
from which we deduce that the projective dimension of $M$ is $\leqslant m$.
But then, replacing the complex 
\[\xymatrix{
0\ar[r] & P^{-k}\ar[r] & P^{-k+1}\ar[r] &\cdots\ar[r] &
P^{\ell-1}\ar[r] & P^{\ell}\ar[r] & 0
}\] 
by a homotopy equivalent one, we may assume $k\leqslant m$. And now the result 
follows immediately.
\eprf

And now we are ready for

\medskip

\nin
\emph{Proof of Lemma~\ref{L1.3}.}\ \ 
Assume $X$ is a quasicompact, quasiseparated scheme. Suppose we are 
given a bounded \tstr\ on
$\dperfs ZX$, and set $\ca={\dperfs ZX}^{\leqslant0}$. If $G$ is a compact
generator, the boundedness of the \tstr\ guarantees that $\T^nG$ belongs to
$\ca$ for some integer $n>0$.
Therefore $G(-\infty,-n]=\{\T^iG\mid i\geqslant n\}$ is contained in  $\ca$,
and we deduce the inclusion
\[
\Dqcs Z(X)_G^{\leqslant-n}\sub\Coprod(\ca)\ .
\]
This much of the proof is pain-free and works without any hypothesis of
noetherianness or finite dimensionality.

Now the category $\dperfs ZX$ is self-dual, with the functor taking
$E\in\dperfs ZX$ to $E^\vee=\RHHom(E,\co_X^{})$ being an equivalence
$\dperfs Z X\cong{\dperfs ZX}\op$. With $\ca={\dperfs ZX}^{\leqslant 0}$ as
in the last paragraph, and setting
$\cb={\dperfs ZX}^{\geqslant 0}$, we have that duality takes the bounded
\tstr\ $(\ca,\cb)$ to the bounded \tstr\ $(\cb^\vee,\ca^\vee)$ on
$\dperfs ZX$. In particular the last paragraph teaches us that
there exists an integer $n>0$ with
\[
\Dqcs Z(X)_G^{\leqslant-n}\sub\Coprod(\cb^\vee)\ .
\]
Theorem~\ref{T27.1}(iii) tells us  
that the standard \tstr\ is in the preferred equivalence
class. Therefore, after increasing $n$ if necessary, we deduce the
inclusion
\[
\Dqcs Z(X)^{\leqslant-n}\sub\Coprod(\cb^\vee)\ ,
\]
where by $\Dqcs Z(X)^{\leqslant-n}$ we mean the standard \tstr. Taking
perpendiculars gives the third inclusion in
\[
\ca^\vee\sub\T(\cb^\vee)^\perp\sub\T\Coprod(\cb^\vee)^\perp\sub
        \Dqcs Z(X)^{\geqslant-n}\ .
\]
And now the assumptions that $X$ is noetherian and finite dimensional
comes in.
If $m$ is the dimension of $X$
then Lemma~\ref{L2.1} tells us
that $\ca\subset\Dqcs Z(X)^{\leqslant m+n}$ and hence
$\Coprod(\ca)\subset\Dqcs Z(X)^{\leqslant m+n}$. This finishes the proof.
\hfill{$\Box$}

\section{Bounded \tstr s on  $\dcohs Z(X)$,
in the presence of dualizing complexes}
\label{S31}

In this section and the next we turn to the problem of bounded 
\tstr s on the category
$\dcohs Z(X)$. By different methods, each section will prove a
theorem showing that (under suitable hypotheses) the bounded
\tstr s on $\dcohs Z(X)$ are all equivalent.

The following lemma
will be useful in both of these sections.

\lem{L31.1}
Let $X$ be a noetherian scheme and let $Z\subset X$ be
a closed subset.
Assume we are given on $\dcohs Z(X)$ a \tstr\ 
$\big(\dcohs Z(X)^{\leqslant0},\dcohs Z(X)^{\geqslant0}\big)$ 
such that 
\[
\dcohs Z(X)=\bigcup_{i=0}^\infty \dcohs Z(X)^{\leqslant i}\ .
\]
Put $\ca=\dcoh(X)^{\leqslant0}$. 

Then $\Coprod(\ca)$ is the aisle of a
compactly generated \tstr.
\elem

For the assertion that $\ca$ generates a \tstr\ on $\Dqcs Z(X)$ 
we could appeal to results in the literature
saying that, if $\ct$ is a good enough triangulated category,
then any small set of objects in $\ct$ generates a \tstr.
The oldest such theorem is Alonso, Jerem{\'\i}as
and Souto~\cite[Proposition~3.2]{Alonso-Jeremias-Souto03},
and the most recent and powerful
is \cite[Theorem~2.3]{Neeman21A}.
But none of these general theorems tells us that
this \tstr\ is compactly generated; and
this compact generation will be important in the next
section.

\prf
Let $G$ be a compact generator of $\Dqcs Z(X)$. Being compact,
$G$ belongs to $\dperfs ZX\subset\dcohs Z(X)$. Because
we are assuming
$\dcohs Z(X)=\cup_{i=0}^\infty \dcohs Z(X)^{\leqslant i}$,
there exists an integer
$n>0$ with $\T^n G\in\dcohs Z(X)^{\leqslant0}=\ca$. Therefore
\[
\Dqcs Z(X)_G^{\leqslant-n}\sub\Coprod(\ca)\ .
\]
Now Theorem~\ref{T27.1}(iii) 
tells us that the standard \tstr\ 
on $\Dqcs Z(X)$ is in the preferred equivalence
class. Therefore, increasing the integer $n$ if
necessary, we deduce 
\be
\item
There exists an integer $n>0$ and a compact generator $G\in\Dqcs Z(X)$ such that
$\T^nG\in\ca$ and
\[
\Dqcs Z(X)^{\leqslant-2n}\sub\Dqcs Z(X)_G^{\leqslant-n}\sub\Coprod(\ca)\ .
\]
\setcounter{enumiv}{\value{enumi}}
\ee

Next: with the integer $n>0$ as in (i), we define 
a full subcategory $\cb\subset\dperfs ZX$ by the formula
\[
\cb\eq\dperfs ZX\cap\Big(\Dqcs Z(X)^{\leqslant-2n}\star\ca\Big)\ .
\]
The suspension functor $\T:\Dqcs Z(X)\la\Dqcs Z(X)$ takes each of
the subcategories
$\dperfs ZX$,  $\Dqcs Z(X)^{\leqslant-2n}$ and $\ca$ into itself,
and the formula defining $\cb$ therefore gives that 
$\T\cb\subset\cb$. Furthermore:
as $\{0\}\subset\Dqcs Z(X)^{\leqslant0}$ we have that
\[
\ca\eq\{0\}\star\ca\sub\Dqcs Z(X)^{\leqslant-2n}\star\ca\ ,
\]
and intersecting with $\dperfs ZX$ we deduce that $\dperfs ZX\cap\ca\subset\cb$.
But (i) gives the first inclusion in
 $G(-\infty,-n]\subset\dperfs ZX\cap\ca\subset\cb$, and
we conclude that $G(-\infty,-n]$ is contained both in $\ca$ and in $\cb$.
Thus 
both $\Coprod(\ca)$ and $\Coprod(\cb)$
must contain 
$\Coprod\big(G(-\infty,-n]\big)=\Dqcs Z(X)_G^{\leqslant-n}$, which contains
$\Dqcs Z(X)^{\leqslant-2n}$. 

From the formula defining $\cb$ we immediately
obtain the first inclusion in
\be
\setcounter{enumi}{\value{enumiv}}
\item
$\qquad\cb\sub\Dqcs Z(X)^{\leqslant-2n}\star\ca\sub\Coprod(\ca)\ ,$
\setcounter{enumiv}{\value{enumi}}
\ee
and the second inclusion is because $\Coprod(\ca)$ contains
both $\Dqcs Z(X)^{\leqslant-2n}$ and $\ca$ and is closed under extensions.

And now we apply Theorem~\ref{T29.1}. 
We remind
the reader: Theorem~\ref{T29.1} says that, for any
object $A\in\dcohs Z(X)$ and for any integer $M$,  there exists an object 
$B\in\dperfs ZX$ and a triangle $D\la B\la A\la\T D$
with $D\in\Dqcs Z(X)^{\leqslant M}$. Putting $M=-2n$ with $n>0$ as in (i), and
applying Theorem~\ref{T29.1} 
to an arbitrary object $A\in\ca\subset\dcohs Z(X)$, the part $D\sr B\sr A$
of the triangle tells us that $B\in\dperfs ZX$ also belongs to
$\Dqcs Z(X)^{\leqslant-2n}\star\ca$, thus $B$ must lie in
$\cb=\dperfs ZX\cap\big(\Dqcs Z(X)^{\leqslant-2n}\star\ca\big)$. But then the part 
$B\sr A\sr\T D$ of the triangle
exhibits $A\in\ca$ as an object
belonging to $\cb\star\Dqcs Z(X)^{\leqslant-2n-1}$. As $A\in\ca$
is arbitrary this gives the first
inclusion in
\be
\setcounter{enumi}{\value{enumiv}}
\item
$\qquad
\ca\sub\cb\star\Dqcs Z(X)^{\leqslant-2n-1}\sub\Coprod(\cb)\ ,$
\setcounter{enumiv}{\value{enumi}}
\ee
and the second inclusion is because $\Coprod(\cb)$ contains
both $\cb$ and $\Dqcs Z(X)^{\leqslant-2n-1}$ and is closed under extensions.

The inclusions (ii) and (iii) combine to tell us that 
$\Coprod(\ca)=\Coprod(\cb)$, and because $\cb$ is
a set of compact objects we have, by Reminder~\ref{R9.5}(i), that
$\Coprod(\ca)=\Coprod(\cb)$ is the aisle of a compactly 
generated \tstr\ on
$\Dqcs Z(X)$.
\eprf

\rmk{R31.999}
If $\big(\dcohs Z(X)^{\leqslant0},\dcohs Z(X)^{\geqslant0}\big)$ is
a bounded \tstr\ on $\dcohs Z(X)$, then the hypothesis 
$\dcohs Z(X)=\cup_{i=0}^\infty \dcohs Z(X)^{\leqslant i}$ is satisfied and
the conclusion of
Lemma~\ref{L31.1} holds. That is: with
$\ca=\dcohs Z(X)^{\leqslant0}$ we have that
$\big(\Coprod(\ca),\Coprod(\ca)^\perp\big)$ is
a compactly generated \tstr.

By Lemma~\ref{L3.1} it follows that, to prove that the 
bounded \tstr s on $\dcohs Z(X)$ are
equivalent, it suffices to show that the \tstr s they generate
on $\Dqcs Z(X)$ are all in the preferred equivalence class.
\ermk

We have proved in this section a preparatory result, 
and it is time to come to the point.
The following result, about bounded \tstr s on
$\dcohs Z(X)$, 
is proved using a minor variant of the technique that
worked to prove Lemma~\ref{L1.3}.

\thm{T31.3}
Let $X$ be a finite-dimensional, noetherian scheme, and assume that $X$ has a
dualizing complex. Let $\big(\dcohs Z(X)^{\leq0},\dcohs Z(X)^{\geq0}\big)$ be any
bounded \tstr\ on $\dcohs Z(X)$.

By Lemma~\ref{L31.1} the set of objects $\dcohs Z(X)^{\leq0}$ generates a 
(compactly generated)
\tstr\ on the category $\Dqcs Z(X)$. We assert that this \tstr\ must be
in the preferred equivalence class. 
\ethm

\prf
Put $\ca=\dcohs Z(X)^{\leq0}$, and 
let $G$ be a compact
generator for $\Dqcs Z(X)$. Then $G\in\dperfs ZX\subset\dcohs Z(X)$,
and as the \tstr\ on $\dcohs Z(X)$ is assumed bounded 
there must exist an integer
$n>0$ such that $\T^n G\in\ca$.
And as in part (i) of the proof of Lemma~\ref{L31.1},
we deduce the inclusion
\[
\Dqcs G(X)_G^{\leqslant-n}\sub\Coprod(\ca)\ .
\]

Now let $D$ be a dualizing complex.
The category $\dcohs Z(X)$ is self-dual, with the functor taking
$E\in\dcohs Z(X)$ to $E^*=\RHHom(E,D)$ being an equivalence
$\dcohs Z(X)\cong{\dcohs Z(X)}\op$. With $\ca={\dcohs Z(X)}^{\leqslant 0}$ as
above, and setting
$\cb={\dcohs Z(X)}^{\geqslant 0}$, we have that duality takes the bounded
\tstr\ $(\ca,\cb)$ to the bounded \tstr\ $(\cb^*,\ca^*)$ on
$\dcohs Z(X)$. In particular the last paragraph teaches us that
there exists an integer $n>0$ with
\[
\Dqcs Z(X)_G^{\leqslant-n}\sub\Coprod(\cb^*)\ .
\]
Theorem~\ref{T27.1}(iii) tells us  
that the standard \tstr\ is in the preferred equivalence
class. Therefore, after increasing $n$ if necessary, we deduce the
inclusion
\[
\Dqcs Z(X)^{\leqslant-n}\sub\Coprod(\cb^*)\ ,
\]
where by $\Dqcs Z(X)^{\leqslant-n}$ we mean the standard \tstr. Taking
perpendiculars gives the third inclusion in
\[
\ca^*\sub\T(\cb^*)^\perp\sub\T\Coprod(\cb^*)^\perp\sub
        \Dqcs Z(X)^{\geqslant-n}\ .
\]
Thus $\ca^*\subset\dcohs Z(X)\cap\Dqcs Z(X)^{\geq-n}=\dcohs Z(X)^{\geq-n}$,
where this time we mean $\dcohs Z(X)^{\geq-n}$ with respect to the
standard \tstr\ on $\dcohs Z(X)$.

And now we recall that the dualizing complex $D$ is a bounded
complex of injectives, and hence there exists an integer $m>0$ such that
the functor $\RHHom(-,D)$ takes $\dcohs Z(X)^{\geqslant-n}$ into
$\dcohs Z(X)^{\leqslant m}$. Thus applying the functor $\RHHom(-,D)$
to the inclusion above delivers that
$\ca\subset\Dqcs Z(X)^{\leqslant m}$, and hence
$\Coprod(\ca)\subset\Dqcs Z(X)^{\leqslant m}$. 

This concludes the proof.
\eprf

\section{Bounded \tstr s on  $\dcoh(X)$, without dualizing complexes}
\label{S3}

We promised two theorems about bounded \tstr s on $\dcohs Z(X)$.
We have so far delivered one: the combination of Theorem~\ref{T31.3} and
Remark~\ref{R31.999} tells us that, in the presence of a
dualizing complex, the bounded \tstr s on $\dcohs Z(X)$ 
are all equivalent. In this section
we will prove the second theorem---the scheme $X$ is no longer assumed to
have a dualizing complex, but the price
we pay is that we can only handle the case $Z=X$, and $X$ must be separated.  

We have often appealed to \cite[Section~1]{Neeman17}, but until 
now we haven't needed the other sections. In particular
we have yet to appeal to the strong generation that is proved in the
article. This will soon change, see the proof of Theorem~\ref{T3.3}
below.
To prepare ourselves it might help to include:

\rmd{R3.5}
Let $\ct$ be a triangulated category with coproducts, and let $\ch\subset\ct$
be a full subcategory. We recall \cite[Definition~1.3(ii)]{Neeman17}:
for integers $n>0$
the full subcategories $\Coprod_n(\ch)$ are defined, inductively, by the
rule
\[
\Coprod_1(\ch)=\Add(\ch),\qquad\Coprod_{n+1}(\ch)=\Coprod_1(\ch)\star\Coprod_n(\ch).
\]
And the unbounded versions of the strong generations theorems
in \cite{Neeman17}, that is \cite[Theorems~2.1 and 2.3]{Neeman17}, are assertions
that, subject to certain hypotheses on the scheme $X$,
there exists an object $H$, in either $\dperfs ZX$ or
in $\dcoh(X)$, and an integer $N>0$ with
\[
\Dqc(X)=\Coprod_N\big(H(-\infty,\infty)\big)
\] 
where $H(-\infty,\infty)$ is as in Reminder~\ref{R9.7}. And 
\cite[Theorem~2.3]{Neeman17} was later improved in
Aoki~\cite[proof of the Main Theorem]{Aoki20},
where the condition on $X$ was 
made less stringent.

Two things should be noted. First of all: the Main Theorem in~\cite{Aoki20}
is a
statement involving only
$\dcoh(X)$, but the proof is by way of the assertion 
above about
$\Dqc(X)$. See \cite[the final paragraph before the bibliography]{Aoki20},
and observe that there is a 
small typo in the display formula of
that paragraph, the right
hand side should be $\Dqc(X)$, not $\dcoh(X)$.
And the careful reader will observe that
the notation in \cite[proof of the Main Theorem]{Aoki20}
is different from that of the older
\cite[Theorem~2.3]{Neeman17}. When we cite the result
the hypothesis on $X$ will be as in 
\cite[proof of the Main Theorem]{Aoki20}, 
but the conclusion
about $\Dqc(X)$ will be expressed in terms of the 
$\Coprod_N\big(H(-\infty,\infty)\big)$
of  \cite[Theorem~2.3]{Neeman17}. After all: we
have put the reader through mountains of notation already and should
spare her unnecessary duplication.
\ermd

The time has now come for Theorem~\ref{T3.3} and its proof.

\thm{T3.3}
Let $X$ be a noetherian, separated, finite-dimensional, quasiexcellent 
scheme. Let $\big(\dcoh(X)^{\leq0},\dcoh(X)^{\geq0}\big)$ be a bounded
\tstr\ on $\dcoh(X)$.

Then the \tstr\ on $\Dqc(X)$, generated by $\dcoh(X)^{\leqslant0}$,
belongs to
the preferred equivalence class.
\ethm

\prf
Put $\ca=\dcoh(X)^{\leqslant0}$, and
choose a compact generator $G\in\Dqc(X)$. As $G$ belongs 
to $\dperf X\subset\dcoh(X)$
and the \tstr\ on $\dcoh(X)$ is bounded, we may choose an integer $n$ with 
$\T^nG\in\ca$. Therefore 
\be
\item
$\qquad \Dqc(X)^{\leq-n}_G\subset\Coprod(\ca)$.
\setcounter{enumiv}{\value{enumi}}
\ee
We need to prove the reverse inclusion.

The hypotheses of the Theorem allow us to 
appeal to 
Aoki~\cite[proof of the Main Theorem]{Aoki20}, which says
that
\be
\setcounter{enumi}{\value{enumiv}}
\item
We may choose and fix an object $H\in\dcoh(X)$ and an integer
$N>0$ with $\Dqc(X)=\Coprod_N\big(H(-\infty,\infty)\big)$;
see Reminder~\ref{R3.5} for an explanation of the notation.
\setcounter{enumiv}{\value{enumi}}
\ee
Because $H$ is an object of $\dcoh(X)$ and the \tstr\ 
$\big(\dcoh(X)^{\leqslant0},\dcoh(X)^{\geqslant0}\big)$ is assumed bounded,
we may choose an integer $m>0$ such that
\[
H[m,\infty)\sub\dcoh(X)^{\geqslant1}\eq\ca^\perp\sub\Coprod(\ca)^\perp\ .
\]
Now the \tstr\ $\big(\Coprod(\ca),\Coprod(\ca)^\perp\big)$
is compactly generated by Lemma~\ref{L31.1}, and
\cite[Proposition~A.2]{Alonso-Jeremias-Souto03}
tells us that
$\Coprod(\ca)^\perp$ is closed
under coproducts (and extensions).
We deduce
\[
\Coprod\big(H[m,\infty)\big)\sub
\Coprod(\ca)^\perp\ .
\]
This gives us the final inclusion in the string
\begin{eqnarray*}
\Dqc(X)^{\geqslant M}&\eq&
\Dqc(X)^{\geqslant M}\cap\Coprod_N\big(H(-\infty,\infty)\big)\\
&\sub&\Coprod\big(H[m,\infty)\big)\\
&\sub&\Coprod(\ca)^\perp\ .
\end{eqnarray*}
By (ii) the equality holds for any integer $M$. The existence of a
sufficiently large integer $M$, for
which  
the middle inclusion holds, comes from
 \cite[Lemma~2.4]{Neeman17}; this integer $M$ depends on $m$, $N$ and $H$.

Taking 
perpendiculars
we deduce
\[
\Coprod(\ca)\sub \Dqc(X)^{\leqslant M-1}\ ,
\]
and together with (i) this concludes the proof.
\eprf

\rmk{R997.3}
Let $X$ be a finite-dimensional noetherian scheme, and 
$Z\subset X$ a closed subset.
In this article we have proved theorems about the 
equivalence classes of  bounded \tstr s on $\dperfs ZX$ and on
$\dcohs Z(X)$. Let us begin with $\dperfs ZX$: in  Lemma~\ref{L1.3}
we proved that the bounded \tstr s on  
$\dperfs ZX$ are all equivalent, unconditionally. Of course
we should read this in conjunction with Theorem~\ref{T0.1},
which tells us that there are no bounded \tstr s on
$\dperfs ZX$ unless $Z$ is contained in the regular
locus of $X$, in which case $\dperfs ZX=\dcohs Z(X)$. 
\ermk

Focusing our attention on the results we have about $\dcohs Z(X)$,
we summarize them as follows.

\rmk{R997.3333098}
The bounded \tstr s on the category $\dcohs Z(X)$ are all equivalent
if any of the following holds:
\be
\item
$Z$ is contained in the regular locus of $X$.
\item
$X$ admits a dualizing complex.
\item
$X$ is separated and quasiexcellent, and the inclusion $Z\subset X$ is an equality.
\ee
For the proof that (i) suffices see the discussion preceeding (\ref{ST897.3.1}).
The proof that (ii) suffices is by combining  Theorem~\ref{T31.3} with 
Remark~\ref{R31.999}.
And the proof that (iii) suffices is by combining Theorem~\ref{T3.3} with
Remark~\ref{R31.999}.
\ermk

Let us note that there is clearly overlap among conditions (i), (ii) and (iii)
of Remark~\ref{R997.3333098}.
For example: assume that the noetherian scheme 
$X$ is J-1, meaning the set $U$ of regular points
in $X$
is open. Then, under the hypothesis in
Remark~\ref{R997.3333098}(i), the closed subset $Z$
must be contained in $U$. But Corollary~\ref{C27.2.1}
delivers an equivalence
$\dcohs Z(U)\cong\dcohs Z(X)$, and $U$ has 
a dualizing complex; after all on a regular scheme the structure sheaf is a
dualizing complex. Thus for J-1 schemes Remark~\ref{R997.3333098}(i) follows
easily from Remark~\ref{R997.3333098}(ii).

Presumably the equivalence of the bounded \tstr s on
$\dcohs Z(X)$ holds in a generality greater than we can prove now.

\section{The standard \tstr\ is in the preferred 
equivalence class}
\label{S27}

We remind the reader of

\medskip

\nin
{\bf Theorem~\ref{T27.1}.}\ \ 
\begin{emph}
Let $X$ be a quasicompact, quasiseparated scheme. Assume $Z\subset X$
is a closed subset such that $X-Z$ is quasicompact. 
Let $\Dqcs{Z}(X)\subset\Dqc(X)$ be as in Reminder~\ref{R9.17}.

Then the following assertions are true:
\be
\item
The compact objects in $\Dqcs{Z}(X)$ are those perfect complexes on
$X$ whose cohomology is supported in $Z$. We will denote the 
full subcategory of such objects $\dperfs{Z}X$.
\item
There exists a single compact object $G\in\Dqcs{Z}(X)$ generating 
$\Dqcs{Z}(X)$.
\item
The standard \tstr\ on $\Dqcs{Z}(X)$ is in the preferred equivalence
class. 
\item
The category $\Dqcs Z(X)$ is weakly approximable, as in
\cite[Definition~0.21]{Neeman17A}.
\ee
\end{emph}

\medskip

In this section we will prove parts (i), (ii) and (iii)
of the theorem; note that parts (i), (ii) and (iii) have
all been
used in the earlier
part of this article, they were crucial to
the proof we have given of Theorem~\ref{T0.1}. In the
next section we prove part (iv) of the theorem, which
is not relevant to the current article but turns out
to be useful in other contexts.

\rmk{R27.2}
The new portion of Theorem~\ref{T27.1} is the assertions
(iii) and (iv). In their current generality, parts (iii)
and (iv) go way beyond what was known.

If $Z=X$ then  parts (i) and (ii) of Theorem~\ref{T27.1} are
both relatively old,
they were proved already in 
Bondal and Van den Bergh~\cite[Theorem~3.1.1]{BondalvandenBergh04}.
The only known result, in the direction of
parts (iii) and (iv) of Theorem~\ref{T27.1}, may
be found in the (recent)~\cite{Neeman17A};
it assumes not only that $Z=X$, but also
restricts $X$ to be separated,

Going back to parts (i) and (ii): more recently
than the old \cite[Theorem~3.1.1]{BondalvandenBergh04},
and now with 
$Z$ allowed to be a proper closed subset of $X$,
assertions (i) and (ii) were proved in
Rouquier~\cite[Theorem~6.8]{Rouquier08}, using essentially
the same technique as in
\cite{BondalvandenBergh04}. A few years later
Hall and Rydh~\cite{Hall-Rydh13} studied
the compact generation of
$\Dqcs{Z}(X)$ when $X$ is an algebraic stack.
If we specialize the results Hall and Rydh~\cite{Hall-Rydh13}
to the case of schemes, then (i) and (ii) are
included---but the method of the
proof is somewhat different.
As we said in the
opening paragraph: (i) and (ii) are not new.

We should say epsilon about what happens when we move
away from schemes.
Theorem~\ref{T27.1}(i) fails in general;   
the compact objects are always perfect,
but for general enough algebraic stacks the converse is false.

And the stacky version of Theorem~\ref{T27.1}(ii) is that,
for quasicompact stacks with quasi-finite and separated diagonal (e.g.\  
algebraic spaces and separated Deligne-Mumford stacks),
the combination of
\cite[Theorem~A and Lemma~4.9]{Hall-Rydh13}
establishes that $\Dqcs Z(X)$ has a single compact generator.

For stacks with infinite stabilizers, the article \cite{Hall-Rydh13}
does contain interesting results about compact generation
in $\Dqcs Z(X)$.
 However, as the reader can
see in Hall and Rydh~\cite{Hall-Rydh15}, there is rarely a single
compact generator for such  $\Dqcs Z(X)$.
\ermk

As we said in Remark~\ref{R27.2}, only parts (iii) and
(iv) of Theorem~\ref{T27.1}
are genuinely new. We will nevertheless give a complete proof, for two reasons:
\be
\item
The relevant beautiful work, in Bondal and
Van den Bergh~\cite{BondalvandenBergh04},
Rouquier~\cite{Rouquier08}
and
in Hall and Rydh~\cite{Hall-Rydh13}, deserves to become
better known than it seems to be.
\item
The exposition we chose is to suitably modify the old proof
of \cite[Theorem~3.1.1]{BondalvandenBergh04}
so that (new) parts (iii) and (iv) of Theorem~\ref{T27.1} come
along for the
ride. Hence the experts might be curious to compare the old proofs
with the new, to see where the modifications
occur. The short summary is that the new proof uses
functorial constructions more systematically, and makes less
use
of the Thomason-Trobaugh localization theorem than the
older approaches. For the reader's convenience, the
part of the
Thomason-Trobaugh theorem which we will use can be
found in Corollary~\ref{C27.3} below.
\ee
We will prove Theorem~\ref{T27.1} through a sequence
of lemmas, attempting to make the argument
reasonably self-contained. Before we start we set

\ntn{N27.7999}
Let $X$ be a quasicompact, quasiseparated scheme, let $Z\subset X$ be a
closed subset with quasicompact complement, and let $U\subset X$ be
a quasicompact open subset. Write $i:U\la X$ for the open immersion,
and put $Z_1=Z-U$.
\entn

We begin with

\lem{L27.2}
With the conventions of Notation~\ref{N27.7999},
the restriction functor $i^*:\Dqcs Z(X)\la\Dqcs{Z\cap U}(U)$
is a Verdier quotient map and has a right adjoint.
Explicitly: the ordinary derived pushforward
map $i_*:\Dqc(U)\la\Dqc(X)$ restricts to the respective
subcategories to deliver a functor which,
by abuse of notation, we also call
$i_*:\Dqcs{Z\cap U}(U)\la\Dqcs{Z}(X)$. This (second) functor
$i_*$ is the requisite right adjoint, it is fully faithful, and
the counit
of the adjunction $\e:i^*i_*\la\id$ is an isomorphism.

Moreover: the kernel of $i^*$ is
the localizing subcategory $\Dqcs{Z_1}(X)\subset\Dqcs Z(X)$.
\elem

\prf
Consider the commutative square of open immersions
\[\xymatrix@C+30pt{
U-Z \ar[r]^-g\ar[d]_f & U\ar[d]^i \\
X-Z \ar[r]_-j & X
}\]
This square is cartesian, and flat base-change gives that
the following square of derived functors commutes
up to natural isomorphism
\[\xymatrix@C+30pt{
\Dqc(U-Z)\ar[d]_{f_*} & \Dqc(U)\ar[d]^{i_*}\ar[l]_-{g^*} \\
\Dqc(X-Z)  & \Dqc(X)\ar[l]^-{j^*}
}\]
Now apply the isomorphic composites to the subcategory 
$\Dqcs{Z\cap U}(U)\subset\Dqc(U)$.
We have
\[
j^*i_*\Dqcs{Z\cap U}(U)\cong f_*g^*\Dqcs{Z\cap U}(U)=0\ ,
\]
where the vanishing is since $g^*$ obviously kills
$\Dqcs{Z\cap U}(U)$. Hence $j^*$ kills $i_*\Dqcs{Z\cap U}(U)$, giving the
second inclusion in
\[
i^*\Dqcs{Z}(X)\subset\Dqcs{Z\cap U}(U)\qquad\text{and}\qquad
i_*\Dqcs{Z\cap U}(U)\subset\Dqcs{Z}(X)\ .
\]
The first inclusion is trivial.

Now: the functor $i^*:\Dqc(X)\la\Dqc(U)$ and its right adjoint
$i_*:\Dqc(U)\la\Dqc(X)$ are explicit, and it
is classical that the counit of the adjunction
$\e:i^*i_*\la\id$ is an isomorphism.
By restricting to the subcategories
we obtain a pair of adjoint functors 
$\xymatrix{i^*\colon\Dqcs{Z}(X)\ar@<0.5ex>[r] &
\Dqcs{Z\cap U}(U)\ar@<0.5ex>[l]\colon i_*
 }$ where the counit of the adjunction $\e:i^*i_*\la\id$ is
still an isomorphism. Now formal facts
about such pairs of adjoints, which the
reader can find (for example) either in~\cite[Chapter~9]{Neeman99}
or  in Krause~\cite{Krause10},
allow us to deduce that the map $i_*$ is fully faithful
while the map $i^*$ is a Verdier quotient map. 

Next we compute the kernel of $i^*$. 
A complex in $\Dqc(X)$ belongs to the subcategory 
$\Dqcs Z(X)$ if its restriction to
$X-Z$ is acyclic. And it belongs to the kernel
of $i^*$ if its restriction to $U$ is also acyclic.
Hence the restriction to $X-Z_1=U\cup(X-Z)$ must be
acyclic. Thus the
kernel of the map
$i^*$  is indeed the localizing subcategory
$\Dqcs{Z_1}(X)\subset\Dqcs{Z}(X)$.
\eprf

The case of Lemma~\ref{L27.2} in which $Z_1$ is empty will come up,
hence we formulate it as a corollary. Now $Z_1=Z-U$ is empty if and only if
$Z\subset U$. And if $Z_1$ is empty then
$\Dqcs{Z_1}(X)$ is the zero category, making the Verdier quotient
map with kernel $\Dqcs{Z_1}(X)$ an equivalence. Thus the corollary becomes

\cor{C27.2.1}
Let $X$ be a quasicompact, quasiseparated scheme, let $Z\subset X$ 
be a closed
subset with quasicompact complement, 
and let $U\subset X$ be a quasicompact open set containing $Z$.
Write $i:U\la X$ for the open immersion.

Then the adjoint maps
$\xymatrix{i^*\colon\Dqcs{Z}(X)\ar@<0.5ex>[r] &
\Dqcs{Z}(U)\ar@<0.5ex>[l]\colon i_*
 }$
are inverse equivalences.
\ecor

\rmk{R27.2.2}
In the situation of Corollary~\ref{C27.2.1}, the functor $i_*$ is 
the usual ``extension
by zero''. To see this observe that, if $\cf\in\Dqcs{Z}(U)$ is any
object, then $i_*\cf\in\Dqcs Z(X)$ satisfies the properties
\be
\item
The counit of adjunction $\e$ provides a (canonical) isomorphism
$i^*i_*\cf\cong\cf$, hence the restriction of $i_*\cf$ to
the open subset $U\subset X$ agrees (canonically) with $\cf$.
\item
$i_*\cf$ vanishes on the open subset $(X-Z)\subset X$. 
After all $i_*\cf$ belongs to
the subcategory $\Dqcs Z(X)\subset\Dqc(X)$ 
of complexes whose restriction to $X-Z$
is acyclic.
\ee
In the situation of Corollary~\ref{C27.2.1} we have that
$X=U\cup(X-Z)$, making (i) and (ii) determine $i_*\cf$.

In particular: if $\cf\in\Dqcs{Z}(U)$ is a perfect complex, 
this permits us to see directly that
$i_*\cf\in\Dqcs Z(X)$ is also a perfect complex. Perfection
is local, from (i) we have that $i_*\cf$ is perfect on $U$,
while (ii) tells us that $i_*\cf$ is perfect on $X-Z$.
\ermk

Lemma~\ref{L27.2} has yet another consequence we will use.

\cor{C27.3}
With the conventions of Notation~\ref{N27.7999},
assume that the assertions of
Theorem~\ref{T27.1}~(i) and (ii)
hold for both categories $\Dqcs{Z_1}(X)$ and $\Dqcs{Z}(X)$.
Consider the restriction functor 
$i^*:\Dqcs{Z}(X)\la\Dqcs{Z\cap U}(U)$. Then, for every compact
object $C\in\Dqcs{Z\cap U}(U)$, there exists a 
compact object
$\wt C\in\Dqcs{Z}(X)$ with $i^*\wt C\cong C\oplus\T C$. 
\ecor

Of 
course: once we have proved Theorem~\ref{T27.1}, the
assumption that Theorem~\ref{T27.1}~(i) and (ii)
hold for both categories $\Dqcs{Z_1}(X)$ and $\Dqcs{Z}(X)$
becomes free.
But we will use the Corollary in the course of the 
proof of Theorem~\ref{T27.1}, using it for
$Z_1\subset Z\subset X$ 
for which we already know that Theorem~\ref{T27.1}~(i) and (ii) apply.

\prf
If we let
\[
\car=\Dqcs{Z_1}(X),\qquad
\cs=\Dqcs{Z}(X)\qquad\text{and}\qquad
\ct=\Dqcs{Z\cap U}(U),
\]
then Lemma~\ref{L27.2} informs us that the functor $i^*:\cs\la\ct$ 
is a Verdier
quotient map with kernel $\car$.

Now: by assumption the assertions of Theorem~\ref{T27.1}~(i) and (ii) hold
for both $\car$ and $\cs$. Thus both categories are compactly
generated (in fact by a single compact object in each case), and
the compact generators in $\car$ are also compact in $\cs$.

The Corollary now follows from
\cite[Corollaries~0.9 and 0.10]{Neeman92A}.
\eprf

\rmk{R27.4}
On a historical note: Corollary~\ref{C27.3} came first, the result in algebraic
geometry preceded the formal facts about triangulated categories that we
used to prove it here. The algebro-geometric statement 
of  Corollary~\ref{C27.3} was
proved by Thomason and Trobaugh~\cite{ThomTro}, and the
original proof was based on the algebro-geometric 
methods developed by Illusie in SGA6.
It came before anyone thought
of applying homotopy-theoretic techniques, such 
as compact objects, in this algebro-geometric
setting.
\ermk

We have recalled the background necessary for the
proof of Theorem~\ref{T27.1}, and the time has come to get 
to work.

\rmk{R27.5}
We begin with an easy part of Theorem~\ref{T27.1}: if $C$ is a 
perfect complex on $X$ which belongs to $\Dqcs{Z}(X)$, then $C$ is
compact in $\Dqcs{Z}(X)$. 

To see this let $C^\vee$ be the dual of the perfect
complex $C$. Then we have an isomorphism 
$\Hom(C,-)\cong H^0(C^\vee\oo-)$, and both
the functor $H^0(-)$ and the functor $C^\vee\oo(-)$ respect coproducts.

We mentioned already the article Hall and Rydh~\cite{Hall-Rydh13},
which points out that this argument breaks down for 
general enough algebraic stacks. The problem
there is that $H^0(-)$ need not commute with coproducts. If $\cs$ is a 
coherent sheaf on $X$ then it is possible for
$H^i(\cs)$ to be nonzero for infinitely many $i$, in which case
the map
\[\xymatrix{
\coprod_{i=0}^\infty H^0\left( \T^i\cs\right)\ar[r] &
H^0\left(\coprod_{i=0}^\infty \T^i\cs\right)
}\]
is not an isomorphism.
\ermk

\lem{L27.7}
Parts (i), (ii) and (iii) of Theorem~\ref{T27.1}
are true whenever $X$ is affine.
\elem

\prf
The proof that $\Dqcs{Z}(X)$
is compactly generated, by a single perfect complex,
may be found
already in \cite[Proposition~6.1]{Bokstedt-Neeman93}. 
We remind the reader: if $X=\spec R$ then
\cite[Step 2 in the proof of Proposition~6.1]{Bokstedt-Neeman93} produces
a single Koszul complex in $\Dqc(X)\cong\D(R)$, namely
\[
B=\bigotimes_{i=1}^r\left(R\stackrel{f_i}\la R\right)\ ,
\]
which we view as sitting in degrees $-r\leq j\leq 0$. And the perfect
complex $B$ is proved to be a compact generator.

Since $B$ is compact and generates $\Dqcs{Z}(X)$, it follows
that all the compact objects in $\Dqcs{Z}(X)$ lie in the thick
subcategory classically generated by
the complex $B$. Thus all
the compacts 
in $\Dqcs{Z}(X)$ must be perfect complexes on $X$.
By Remark~\ref{R27.5}
all the perfects are compact, and we have 
proved (i) and (ii) for affine $X$.
It remains to prove (iii).

Now $B$ is a single compact generator 
and belongs to $\Dqcs{Z}(X)^{\leq0}$ for the standard \tstr. 
Hence 
$B(-\infty,0]\subset\Dqcs{Z}(X)^{\leq0}$ and therefore
$\Dqcs{Z}(X)^{\leq0}_B\subset\Dqcs{Z}(X)^{\leq0}$.
We need to prove the reverse inclusion. 

Take any object $A\in\Dqcs{Z}(X)^{\leq0}$. Because 
$\big(\Dqcs{Z}(X)^{\leq0}_B,\Dqcs{Z}(X)^{\geq0}_B\big)$
is a \tstr\ on $\Dqcs{Z}(X)$, we can form a triangle
$A'\la A\la A''\la \T A'$ with $A'\in
\Dqcs{Z}(X)^{\leq0}_B\subset\Dqcs{Z}(X)^{\leq0}$
and with $A''\in\Dqcs{Z}(X)^{\geq1}_B$. The portion of the
triangle
$A\la A''\la \T A'$  expresses $A''$ as an extension 
with both $A$ and $\T  A'$ belonging to
$\Dqc(X)^{\leq0}$. Hence
$A''$ must belong to $\Dqcs{Z}(X)^{\leq0}\cap\Dqcs{Z}(X)^{\geq1}_B$.
And it clearly suffices to show that $A''=0$.

Assume therefore
that $A$ is a nonzero
object in $\Dqcs{Z}(X)^{\leq0}\cap\Dqcs{Z}(X)^{\geq1}_B$, 
and we will prove a contradiction.
Using the equivalence of categories $\Dqc(X)\cong\D(R)$, we
have that $A$ belongs to $\D(R)^{\leq0}$ and has its cohomology 
supported in $Z\subset X$. Since $A\neq0$ we must have that $H^{-i}(A)\neq0$
for some $i\geq0$. Therefore there exists, in $\Dqc(X)\cong\D(R)$,
a nonzero map $\T^{i} R\la A$. But $A$ has its cohomology supported 
in $Z$, and \cite[Lemma~5.7 coupled with Example~5.6]{Lipman-Neeman07}
tells us that there exists an object $B_n\in\COprod\big(B[0,0]\big)$ 
such that
the nonzero map $\T^{i}R\la A$ factors as $\T^{i}R\la\T^{i}B_n\la A$.
The interested reader can find a review
of the Koszul complexes $B_n$
in Construction~\ref{C28.2}(ii),
and in the concluding paragraph of
Construction~\ref{C28.2}(ii) we recall the proof
that
$B_n\in\COprod\big(B[0,0]\big)$.
Anyway:
This produces for us a nonzero map $\T^{i}B_n\la A$
with
$\T^{i}B_n\in \Dqcs{Z}(X)^{\leq0}_B$,
contradicting
the hypothesis that $A\in\Dqcs{Z}(X)^{\geq1}_B$.
\eprf

In the remainder of this section, and in the next, we will
find ourselves in the following setup.

\con{C27.505}
Let $X$ be a quasicompact, quasiseparated scheme. Let $U$, $V$
be quasicompact open subsets with $X=U\cup V$, and let $Z\subset X$
be a closed subset with quasicompact complement.
\be
\item
Define
\[
W=U\cap V,\qquad Z_U=Z\cap U,\qquad Z_V=Z\cap V,\quad Z_W=Z\cap W.
\]
\item
Let $i:U\la X$, $j:W\la X$ and $k:V\la X$ be the open immersions.
\item
We also want a name for $Z_1=Z_U-W$. Note that
\[
Z_1\eq Z_U-W\eq (Z\cap U)- (V\cap U)\eq (Z-V)\cap U\eq Z-V
\]
\ee
where the last equality is since $Z-V\subset U$.
\econ

The next Lemma is also originally due to
Thomason and Trobaugh~\cite{ThomTro}.

\lem{L27.987}
In the situation of Construction~\ref{C27.505}, assume that
Theorem~\ref{T27.1}~(i) and (ii) hold for the categories
$\Dqcs {Z_1}(U)$ and $\Dqcs{Z_U}(U)$.
Let $P\in\Dqcs{Z_V}(V)$ be a perfect complex.

Then there exists a perfect complex $Q\in\Dqcs Z(X)$ 
and an isomorphism $k^*Q\cong P\oplus\T P$.
\elem

\prf
We are given a perfect
complex $P\in\Dqcs{Z_V}(V)$, and if we
restrict it to $W\subset V$ we obtain a perfect complex
$\wt P\in\Dqcs{Z_W}(W)$. By hypothesis, parts (i) and (ii) of 
Theorem~\ref{T27.1} hold for
$\Dqcs{Z_1}(U)$ and $\Dqcs{Z_U}(U)$, and hence
we may apply
Corollary~\ref{C27.3} to the scheme $U$ with open subset $W$ and
closed subset $Z_U$. 
The perfect complex
$\wt P\oplus\T\wt P\in\Dqcs{Z_W}(W)$
can be expressed as the restriction to $W$ of
a perfect complex $\wt Q\in\Dqcs{Z_U}(U)$.
We have a perfect complex $\wt Q$ on $U$, a perfect complex
$P\oplus\T P$ on $V$, and an isomorphism on $W$ of the two
restrictions; both are isomorphic to $\wt P\oplus\T\wt P$.

We can glue, meaning that in $\Dqcs{Z}(X)$ we may form the triangle
\[\xymatrix{
Q\ar[r] & i_*\wt Q\oplus k_*(P\oplus\T P)\ar[r] & 
j_*(\wt P\oplus\T\wt P)\ar[r] &\T Q\ .
}\]
The restriction of $Q$ to $U\subset X$ is isomorphic to the perfect
complex $\wt Q$, while the restriction of $Q$ to $V\subset X$ is 
isomorphic to the perfect complex $P\oplus\T P$. 
We  deduce that $Q$ is a perfect complex belonging to
$\Dqcs{Z}(X)$, and satisfies $k^*Q\cong P\oplus\T P$.
\eprf

\medskip

\nin
\emph{Proof of parts (i), (ii) and (iii) of Theorem~\ref{T27.1}.}\ \ 
Since $X$ is quasicompact it can be covered by finitely many affine
open subsets, and the proof
will be by induction on the number
$n$ of open sets needed to cover. Lemma~\ref{L27.7} starts the
induction. Thus we may assume
we know that the conclusion of Theorem~\ref{T27.1}
holds for all $X$ such that $n=1$.

Assume that $n>0$ is
given and that parts (i), (ii) and (iii) of the theorem have been
proved for any quasicompact,
quasiseparated scheme $X$ admitting a cover by 
$\leq n$ affine open
subsets. Choose any 
$X$ admitting a cover by $n+1$ open affines.
We need to prove that
parts (i), (ii) and (iii)
of the theorem also hold for X. 

We may write $X$ as 
$X=U\cup V$, with $U$ an affine open subset and $V$ an open subset
admitting a 
cover by $n$ affine open subsets. Let $Z\subset X$ be a closed 
subset with quasicompact complement.
By induction parts (i) and (ii) of the theorem
hold for the inclusion of the closed
subset $Z_V\subset V$, with $Z_V=Z\cap V$ as in
Construction~\ref{C27.505}. We may therefore choose a compact
generator $P\in\Dqcs{Z_V}(V)$ which is a perfect complex
on $V$. Now $U$ is affine,
and parts (i) and (ii) of 
Theorem~\ref{T27.1} have already been proved for
the categories $\Dqcs{Z_1}(U)$ and $\Dqcs{Z_U}(U)$.
Lemma~\ref{L27.987} allows us to find a
perfect complex $Q\in\Dqcs Z(X)$ with
$k^*Q\cong P\oplus\T P$.
In particular $k^*Q$ is a compact generator
for $\Dqcs{Z_V}(V)$.

We have already observed that
Theorem~\ref{T27.1}~(i) and (ii)
hold for the inclusion
$Z_1\subset U$, permitting us to choose a single perfect complex
$H\in\Dqcs{Z_1}(U)$ which generates. The inclusion $i:U\la X$ allows
us to define $S=i_*H$,
which by Remark~\ref{R27.2.2} is a perfect
complex on $X$ belonging to $\Dqcs{Z_1}(X)$.
We define $G=Q\oplus S$. The object $G$ is a perfect
complex belonging to $\Dqcs Z(X)$, and we will
now analyze $G$ to prove
Theorem~\ref{T27.1}~(i), (ii) and (iii)
for the inclusion $Z\subset X$. We will show that
$G$ is a compact generator for $\Dqcs{Z}(X)$,
that the standard
\tstr\ on $\Dqcs{Z}(X)$ is equivalent to the \tstr\ 
$\big(\Dqcs{Z}(X)^{\leq0}_G,\Dqcs{Z}(X)^{\geq0}_G\big)$, and that
the compact objects in $\Dqcs Z(X)$ coincide with the perfect
complexes.

The trivial part is the observation that, since $G$ is perfect, it
is compact and belongs to $\Dqcs{Z}(X)^{\leq m}$ for the
standard \tstr\ 
and for some $m>0$. Hence 
$\Dqcs{Z}(X)^{\leq 0}_G\subset\Dqcs{Z}(X)^{\leq m}$. In the coming
paragraphs we
proceed to prove that, possibly after increasing the integer $m$,
we can guarantee the
non-trivial inclusions
\[
\Dqcs{Z}(X)^{\leq 0}\subset \Coprod\big(G(-\infty,m]\big)\qquad\text{and}\qquad
\Dqcs{Z}(X)\subset\Coprod\big(G(-\infty,\infty)\big)\ .
\]
  
Now recall that $G=Q\oplus S$
with $S$ supported on $Z_1=Z-V$. Hence $k^*$ annihilates
$S$, and we have
\[ 
k^*G\cong k^*(Q\oplus S)
\cong k^* Q
\]
where by construction we know that $k^*Q=k^*G$ is a compact
generator for $\Dqcs{Z_V}(V)$.
We deduce
\begin{enumerate}
\item
$\qquad\Dqcs{Z_V}(V)=\Coprod\big(k^*G(-\infty,\infty)\big)\ .$
\setcounter{enumiv}{\value{enumi}}
\end{enumerate}
Because Theorem~\ref{T27.1}(iii) holds
for $Z_V\subset V$ we may assume, after shifting $G$, that 
\begin{enumerate}
\setcounter{enumi}{\value{enumiv}}
\item
$\qquad\Dqcs{Z_V}(V)^{\leq0}\subset\Coprod\big(k^*G(-\infty,0]\big)\ .$
\setcounter{enumiv}{\value{enumi}}
\end{enumerate}

Now $G=Q\oplus S$ belongs to
$\Coprod\big(G(-\infty,0]\big)$, which by Remark~\ref{R9.4} is closed under direct
summands. Hence $S$ belongs to 
$\Coprod\big(G(-\infty,0]\big)\subset\Coprod\big(G(-\infty,\infty)\big)$.
But $S\in\Dqcs{Z_1}(X)\cong\Dqcs{Z_1}(U)$ was chosen to be a
compact generator, and Theorem~\ref{T27.1}(iii) holds also
for the inclusion $Z_1\subset U$. Hence we may, after shifting $G$ 
some more,
assume further the first inclusion below 
\begin{enumerate}
\setcounter{enumi}{\value{enumiv}}
\item
$\qquad\Dqcs{Z_1}(X)^{\leq0}\subset\Coprod\big(S(-\infty,0]\big)
\subset\Coprod\big(G(-\infty,0]\big)\ .$ 
\setcounter{enumiv}{\value{enumi}}
\end{enumerate}
The second inclusion
is because  
$S$ belongs $\Coprod\big(G(-\infty,0]\big)$, see above.

Similarly we obtain inclusions
\begin{enumerate}
\setcounter{enumi}{\value{enumiv}}
\item
$\qquad      
\Dqcs{Z_1}(X)\subset\Coprod\big(S(-\infty,\infty)\big)
\subset\Coprod\big(G(-\infty,\infty)\big)\ .
$
\setcounter{enumiv}{\value{enumi}}
\end{enumerate}
Finally we choose an integer $N>0$ such that
\begin{enumerate}
\setcounter{enumi}{\value{enumiv}}
\item
$G$ belongs to $\Dqcs{Z}(X)^{\leq N}$, and 
$k_*\Dqcs{Z_V}(V)^{\leq0}\subset\Dqcs{Z}(X)^{\leq N}$.
\setcounter{enumiv}{\value{enumi}}
\end{enumerate}

We are done making choices: we have constructed the perfect complex 
$G\in\Dqcs Z(X)$ satisfying properties (1) through (5). The strategy
of the remainder of the proof will first be to make
estimates, telling us about the containment
in categories of the form $\genu G{}{A,B}$ of objects
of the form $k_*k^*M$. We will then combine
these with triangles of the form $L\la M\la k_*k^*M$, and with
the estimates we already have in (3) and (4), to inform us
for which $A<B$ a category $\genu G{}{A,B}$ can be
expected to contain an object $M\in\Dqcs Z(X)$.

We introduced enough strategy, it is time to get to work.
Complete the unit of adjunction $G\la k_*k^*G$ to a triangle
$L\la G\la k_*k^*G\la\T L$ in the category
$\Dqcs{Z}(X)$. By (5) the
object $G$ belongs
to $\Dqcs{Z}(X)^{\leq N}$ and hence $k^*G$ belongs to 
$\Dqcs{Z_V}(V)^{\leq N}$. Applying (5) again
we deduce that $k_*k^*G$ belongs to $\Dqcs{Z}(X)^{\leq2N}$.
The triangle $k_*k^*G\la\T L\la \T G$ expresses $\T L$ as an
extension of two objects in $\Dqc(X)^{\leq2N}$, hence $\T L$ must
lie in $\Dqc(X)^{\leq2N}$. But $\T L$ is annihilated by $k^*$,
and Lemma~\ref{L27.2} tells us that the kernel of
$k^*$ is $\Dqcs{Z_1}(X)$. Thus the crude estimates above
combine to give
\begin{enumerate}
\setcounter{enumi}{\value{enumiv}}
\item
In the triangle $L\la G\la k_*k^*G\la\T L$ we have
\[
\T L\in\Dqcs{Z_1}(X)^{\leq2N}\subset\Coprod\big(G(-\infty,2N]\big)
\]
\setcounter{enumiv}{\value{enumi}}
\end{enumerate}
where the inclusion is by (3) above.

And now the triangle 
$G\la k_*k^*G\la\T L$ exhibits $k_*k^*G$ as an extension of two
objects, which by (6) both belong to
$\Coprod\big(G(-\infty,2N]\big)$. Since $\Coprod\big(G(-\infty,2N]\big)$
    is closed under extensions,
it follows that
\begin{enumerate}
\setcounter{enumi}{\value{enumiv}}
\item
The object $k_*k^*G$ belongs to
$\Coprod\big(G(-\infty,2N]\big)$.
\setcounter{enumiv}{\value{enumi}}
\end{enumerate}
The above combines to give
\begin{enumerate}
\setcounter{enumi}{\value{enumiv}}
\item
We have the string of inclusions
\begin{eqnarray*}
k_*k^*\Dqcs{Z}(X)^{\leq0}&\subset& k_*\Dqcs{Z_V}(V)^{\leq0}\\
 &\subset&
k_*\Coprod\big(k^*G(-\infty,0]\big)\\
 &\subset&
\Coprod\big(k_*k^*G(-\infty,0]\big)\\
&\subset&\Coprod\big(G(-\infty,2N]\big)\ ,
\end{eqnarray*}
\setcounter{enumiv}{\value{enumi}}
\end{enumerate}
where the second inclusion is by
(2), the third inclusion is 
because $k_*$ respects coproducts and extensions, and
the last inclusion is by (7).

The inclusions below are similar.
\begin{enumerate}
\setcounter{enumi}{\value{enumiv}}
\item
We leave to the reader the string of inclusions
\begin{eqnarray*}
k_*k^*\Dqcs{Z}(X)&\subset& k_*\Dqcs{Z_V}(V)\\
 &\subset&
k_*\Coprod\big(k^*G(-\infty,\infty)\big)\\
 &\subset&
\Coprod\big(k_*k^*G(-\infty,\infty)\big)\\
&\subset&\Coprod\big(G(-\infty,\infty)\big)\ .
\end{eqnarray*}
\setcounter{enumiv}{\value{enumi}}
\end{enumerate}

Now let $A\in\Dqcs{Z}(X)$ be arbitrary, and complete
the unit of adjunction $A\la k_*k^*A$
to the distinguished triangle
$B\la A\la k_*k^*A\la\T B$. Then $B$ is annihilated by $k^*$ and
therefore belongs to its kernel 
$\Dqcs{Z_1}(X)\subset\Coprod\big(G(-\infty,\infty)\big)$,
where the inclusion is by (4). And $k_*k^*A$ belongs to
$k_*k^*\Dqcs{Z}(X)\subset\Coprod\big(G(-\infty,\infty)\big)$, where this time
the inclusion is by (9). As $\Coprod\big(G(-\infty,\infty)\big)$ is closed
under extensions we have that $A\in\Coprod\big(G(-\infty,\infty)\big)$. As
$A\in\Dqcs{Z}(X)$ is arbitrary we deduce that the perfect complex $G$ generates
the category $\Dqcs Z(X)$. In Remark~\ref{R27.5} we noted
that every perfect complex is compact, in
particular $G$ is compact and generates $\Dqcs Z(X)$.
It therefore also classically generates
the thick subcategory of compact objects, proving that every compact object
is a perfect complex. This completes the proof of Theorem~\ref{T27.1}~(i) and (ii).

Now assume $A\in\Dqcs{Z}(X)^{\leq0}$. Then 
$k_*k^*A\in\Coprod\big(G(-\infty,2N]\big)$
by (8). And as $k^*A\in\Dqcs{Z_V}(V)^{\leq0}$ we obtain from (5) that
$k_*k^*A\in\Dqcs{Z}(X)^{\leq N}$, and the triangle 
$k_*k^*A\la\T B\la \T A$
says that $\T B$ must belong to 
$\Dqcs{Z_1}(X)^{\leq N}\subset\Coprod\big(G(-\infty,N]\big)$, where the inclusion is
by (3). In the triangle $B\la A\la k_*k^*A$ we now have that
$B\in\Coprod\big(G(-\infty,N+1]\big)\subset\Coprod\big(G(-\infty,2N]\big)$
and $k_*k^*A$ also belongs to $\Coprod\big(G(-\infty,2N]\big)$. Hence
$A\in\Coprod\big(G(-\infty,2N]\big)$, and as 
$A\in\Dqcs{Z}(X)^{\leq0}$ is arbitrary
we deduce the inclusion $\Dqcs{Z}(X)^{\leq0}\subset\Coprod\big(G(-\infty,2N]\big)$. This concludes the proof of Theorem~\ref{T27.1}(iii).
\hfill{$\Box$}

\section{Proof of the weak
approximability of $\Dqcs Z(X)$}
\label{S28}

The final two sections, the present one and the next,
are written for the experts in approximable triangulated
categories, and can safely be skipped by anyone else.
The proof of Theorem~\ref{T0.1} (about bounded \tstr s
on $\dperfs ZX$), as well as the proofs of
Theorems~\ref{T31.3} and \ref{T3.3} (about bounded \tstr s on
$\dcohs Z(X)$), do not depend
on anything that is yet to come.

\rmk{R27.1}
We promised the reader that, in this section, we will prove part (iv) of
Theorem~\ref{T27.1}. This means that we need to prove the category
$\Dqcs Z(X)$ to be weakly approximable.

By \cite[Definition~0.21]{Neeman17A} this is the assertion
that there exist, in $\Dqcs Z(X)$, a compact object $G$, an integer
$A>0$ and a
\tstr\ satisfying some properties. As the reader will see, the \tstr\
we will use in the proof is the standard \tstr\ on
$\Dqcs Z(X)$.

Combining this with  \cite[Proposition~2.4]{Neeman17A}
allows us to deduce that the standard \tstr\ must be in the
preferred equivalence class. Thus the proof we are embarking on,
of Theorem~\ref{T27.1}(iv), will provide a
second proof of Theorem~\ref{T27.1}(iii).
As the reader will soon see,
the proof of Theorem~\ref{T27.1}(iv)
amounts to a more careful
analysis of the various inclusions that came up
in the proof of Theorem~\ref{T27.1}~(i), (ii) and (iii).

Finally note that, by Theorem~\ref{T27.1}(ii), the category
$\ct=\Dqcs Z(X)$
has a
single compact generator $G$. By Theorem~\ref{T27.1}(i) this compact
generator is a perfect complex. Hence, for the standard \tstr\ $\tst\ct$,
we may choose an integer $A>0$ with $G\in\ct^{\leq A}$ and with
$\Hom(G,\ct^{\leq-A})=0$. What needs proof is the approximability
statement: we need to show that, after
increasing the integer $A$ if necessary,
every object $F\in\ct^{\leq0}$ admits
a triangle $E\la F\la D$ with $E\in\ogenu G{}{-A,A}$ and with
$D\in\ct^{\leq-1}$.
\ermk

\rmk{R28.1.5}
Since the reader of this section is assumed to have some familiarity
with the previous work of the author, we will freely pass
between the categories $\Coprod\big(G[-M,M]\big)$ and the categories
$\ogenu G{}{-M,M}$ of
\cite[Reminder~0.8(xi)]{Neeman17A}.
From \cite[Lemma~1.10]{Neeman17} we know that,
up to changing the integer $M$ a little, these categories
are interchangeable.
\ermk

\con{C28.2}
Let $U=\spec R$ be an affine scheme, and let $Z_1\subset U$ be a
closed subset with quasicompact complement. The quasicompact
open subset $U-Z_1$ can be written as a finite union of open subsets
$\spec{R\left[\frac1{f_i}\right]}$. Choose such a collection $\{f_1,f_2,\ldots,f_r\}$.
\be
\item
For any integer $n>0$,
we define the complex $A_{i,n}$ to be the
mapping cone of the morphism $f_i^n:R\la R$. That is $A_{i,n}$
is the cochain
complex
\[\xymatrix{
0\ar[r] & R\ar[r]^-{f_i^n} & R\ar[r] & 0
}\]
concentrated in degrees $\{-1,0\}$.
\item
We let $B_n$ be the Koszul complex
\[
B_n\eq \bigotimes_{i=1}^r A_{i,n}\ .
\]
In passing we remind ourselves that,
by \cite[Proposition~6.1]{Bokstedt-Neeman93}, the
object $B=B_1$
is a compact generator for $\Dqcs {Z_1}(U)$.
We have already had occasion to appeal to this,
back in the proof
of Lemma~\ref{L27.7}.
\ee
Now the triangles $A_{i,1}\la A_{i,n+1}\la A_{i,n}\la$ exhibit
the containment $\{A_{i,n+1}\}\in \{A_{i,1}\}\star \{A_{i,n}\}$,
and induction tells us that, for all $n$, we must have
$A_{i,n}\in\COprod\big(A_{i,1}[0,0]\big)$.
Tensoring over $1\leq i\leq r$,
we deduce that
\[
B_n\in\COprod\big(B_1[0,0]\big)\ .
\]
\econ

\lem{L28.999}
Let $U$ be an affine scheme, and let $Z_1\subset U$ be a closed
subset with quasicompact complement. Write $(U-Z_1)$ for the
complement, and let $j:(U-Z_1)\la U$ be the
open immersion. In the category $\Dqc(U)$ we can complete
the unit of adjunction $\co_U^{}\la j_*j^*\co_U^{}$ to
a distinguished
triangle $L'\la \co_U^{}\la j_*j^*\co_U^{}\la\T L'$.

Choose any compact generator $H\in\Dqcs{Z_1}(U)$.  With $L'$
as above, we assert that
there exists an integer $r>0$ such that
$L'$ belongs the subcategory $\ogenu H{}{-r,r}\subset\Dqcs{Z_1}(U)$.
\elem

\prf
For the proof of this Lemma we adopt
the notation of Construction~\ref{C28.2}:
in particular $U=\spec R$,
and in the ring $R$ we
are given a finite set of elements $\{f_1,f_2,\ldots,f_r\}$
as in Construction~\ref{C28.2}.

The equivalence of
categories $\Dqc(U)\cong\D(R)$ takes the triangle
$L'\la \co_U^{}\la j_*j^*\co_U^{}\la\T L'$ to a well-understood
triangle $L'\la R\la\ph_*\ph^* R\la\T L'$. In
the notation of Construction~\ref{C28.2},
there is an isomorphism
\[
L'\cong\bigotimes_{i=1}^r\left(R\la R\left[\frac1{f_i}\right]\right)\ .
\]
where the complex $R\la R\left[\frac1{f_i}\right]$ is concentrated
in degrees $\{0,1\}$.
Thus, still in the notation of Construction~\ref{C28.2},
$L'$ is the homotopy colimit of the complexes
\[
\T^{-r}B_n=\bigotimes_{i=1}^r \left(\Tm A_{i,n}\right)\ ,
\]
with $A_{i,n}$ being the mapping cone on the map
$f^n_i:R\la R$ as in Construction~\ref{C28.2}.
Therefore there exists a distinguished triangle
\[\xymatrix{
\ds\coprod_{n=1}^\infty \T^{-r}B_n\ar[r] &
\ds\coprod_{n=1}^\infty \T^{-r}B_n\ar[r] &
L'\ar[r] &  
\ds\coprod_{n=1}^\infty \T^{-r+1} B_n\ .
}\]
From Construction~\ref{C28.2} we recall that, for each $n>0$, we have
$B_n\in\genu{B_1}{}{0,0}=\genu B{}{0,0}$. The triangle therefore
exihibits $L'$ as an object in
$\ogenu B{}{r,r}\star\ogenu B{}{r-1,r-1}\subset\ogenu B{}{-r,r}$.

We have proved the existence of the integer $r>0$, asserted in the
Lemma, for the specific
compact generator $B\in\Dqcs{Z_1}(U)$ of Construction~\ref{C28.2}.
But, up to changing the integer $r>0$, this proves the existence for
any compact generator $H\in\Dqcs{Z_1}(U)$.
\eprf

Because the next Lemma has turned out to be useful in contexts 
only tangentially related to this article, we state and prove 
it in the generality in which it is easy to cite and apply.
The conventions are as in Notation~\ref{N27.7999}, but
because this Lemma has been cited elsewhere we give
a self-contained statement.

\lem{L28.1000}
Let $X$ be a quasicompact, quasiseparated scheme, and let $Z\subset X$
be a closed subset with quasicompact complement. Suppose $V\subset X$
is a quasicompact open subset, and let $k:V\la X$ be the open immersion.
Put $Z_1=Z-V$.

We prove the assertions below.

For any compact generators $G\in\Dqcs Z(X)$ and $H\in\Dqcs{Z_1}(X)$---whose
existence is guaranteed because we already know Theorem~\ref{T27.1}(ii)---the
following holds:
\be
\item
For any compact object $F\in\Dqcs{Z\cap V}(V)$, there exists an integer
$M>0$ with $k_*F\in\ogenu G{}{-M,M}$.  
\item
Take any compact object $F\in\Dqcs Z(X)$, and complete 
the unit of adjunction $F\la k_*k^*F$ to a triangle
$L\la F\la k_*k^*F\la\T F$. Then, for
some integer $M>0$, the object $L\in\Dqcs{Z_1}(X)$ belongs
to the subcategory $\ogenu H{}{-M,M}\subset\Dqcs{Z_1}(X)$.
\setcounter{enumiv}{\value{enumi}}
\ee
\elem

\prf
The first reduction is the following.
\be
\setcounter{enumi}{\value{enumiv}}
\item
Assume we can write the
open immersion
$k:V\la X$ as the composite of two open immersions
$V\stackrel\ell\la W\stackrel m\la X$, and further assume that,
for each of the open immersions $\ell$ and $m$, (i) and (ii) hold
for all closed subsets $Z$ with quasicompact complements. Then 
i) and (ii) hold for
the composite $k=m\ell$, again with any choice of closed subsets.
\ee

Let us prove (iii). We begin with (i) for $k=m\ell$,
assuming that (i) holds for each of $m$ and $\ell$.

Choose therefore
any compact object $F\in\Dqcs{Z\cap V}(V)$, and any compact generator
$\wt G\in\Dqcs{Z\cap W}(W)$. By (i) for the open immersion $\ell$
we have that there exists an integer $M>0$ with
$\ell_*F\in\ogenu{\wt G}{}{-M,M}$. But now we are
given a compact
generator $G\in\Dqcs Z(X)$ and by (i), applied to the open
immersion $m:W\la X$, we have that there exists an integer $N>0$ with
$m_*\wt G\in\ogenu G{}{-N,N}$. Combining,
we deduce that
\[
k_*F=m_*\ell_*F\in m_*\ogenu{\wt G}{}{-M,M}\subset\ogenu{m_*\wt G}{}{-M,M}
\subset\ogenu G{}{-M-N,M+N}\ .
\]

To complete the proof of (iii) we need to show
that (ii) holds for $k=m\ell$, assuming that both (i) and (ii) hold
for each of $k$ and $\ell$.

Let $F$ be a compact object in $\Dqcs Z(X)$. Then $m^*F$ is
a compact object of $\Dqcs{Z\cap W}(W)$, and (ii) applies to
the two units of adjunction $F\la m_*m^*F$ and $m^*F\la\ell_*\ell^*m^*F$. 
Each of the
units of adjunction can be completed to a distinguished
triangle, the first in the category $\Dqcs Z(X)$ and the
second in the category $\Dqcs{Z\cap W}(W)$, which we write as
\[
L_2\la F\la m_*m^* F\la\T L_2,\qquad\qquad
L_3\la m^*F\la \ell_*\ell^*m^*F\la \T L_3\ .
\]
Let $Z_2=Z-W$ and let
\[Z_3=(Z\cap W)-V=(Z-V)\cap W=Z_1\cap W\ .
\]
Now let $H_2\in\Dqcs{Z_2}(X)$ and $H_3\in\Dqcs{Z_3}(W)$ be compact
generators and
the hypothesis (ii), for the open
immersions $\ell$ and $m$,
gives us that we may choose an integer $M>0$ with
$L_2\in\ogenu{H_2}{}{-M,M}$ and with $L_3\in\ogenu{H_3}{}{-M,M}$.
But (i), applied to the open immersion $m:W\la X$ and
the closed subset $Z_1\subset X$, tells us that, possibly after
increasing the integer $M>0$, we can also guarantee that
$m_*H_3\in\ogenu H{}{-M,M}$, where $H\in\Dqcs{Z_1}(X)$ is the compact
generator as in the statement of the Lemma. And as $H\in\Dqcs{Z_1}(X)$
is a compact generator and $H_2\in\Dqcs{Z_2}(X)\subset\Dqcs{Z_1}(X)$
is a perfect complex, the integer $M>0$ can be increased to guarantee
that $H_2\in\genu H{}{-M,M}$. This combines to give
\[
\begin{array}{ccccc}
L_2&\in&\ogenu{H_2}{}{-M,M}&\subset&\ogenu{H}{}{-2M,2M}\ ,\\
m_*L_3&\in& m_*\ogenu{H_3}{}{-M,M}&\subset&\ogenu{H}{}{-2M,2M}\ .
\end{array}
\]
Now form the octahedron
\[\xymatrix{
F\ar@{=}[d]\ar[r] & m_*m^* F \ar[r]\ar[d] & \T L_2\ar[d] \\
F\ar[r] & m_*\ell_*\ell^*m^*F\ar[r]\ar[d] & \T L\ar[d] \\
  & \T m_*L_3\ar@{=}[r] & \T m_*L_3
}\]
The unit of adjunction $F\la k_*k^*F$ is identified with the composite
$F\la m_*m^*F\la m_*\ell_*\ell^*m^* F$ of the octahedron above,
the middle row
of the octahedron provides a triangle
$L\la F\la k_*k^*F\la\T L$, and the
third column provides a triangle $L_2\la L\la m_*L_3$, which
exihibits $L$ as belonging to
$\ogenu{H}{}{-2M,2M}*\ogenu{H}{}{-2M,2M}=\ogenu{H}{}{-2M,2M}$.

With (iii) established, note that we may assume that $X=U\cup V$
with $U\subset X$ an open affine---after all, the quasicompactness of $X$
guarantees that it is a
finite union of affine open subsets, and (iii) allows
us to add them one at a time.

We begin with the proof of (ii).
Let $F$ be any perfect complex in $\Dqcs Z(X)$, and complete
the unit of adjunction $\eta:F\la k_*k^*F$ to a triangle
$L\la F\la k_*k^*F\la\T L$. 
Now $L$ is in the kernel of the functor
$k^*:\Dqcs Z(X)\la\Dqcs{Z\cap V}(V)$, and Lemma~\ref{L27.2}
informs us that this kernel is $\Dqcs{Z_1}(X)$.
And the object $H\in\Dqcs{Z_1}(X)$ was chosen, in the
statement of the Lemma, to be a compact generator.
By
Corollary~\ref{C27.2.1} the restriction functor
$i^*:\Dqcs{Z_1}(X)\la\Dqcs{Z_1}(U)$ is an equivalence of
categories. Assertion (ii) of the Lemma, which we want to
prove, is that there exists an integer $M>0$ with
$L\in\ogenu H{}{-M,M}$, and by applying the equivalence
$i^*$ 
it suffices to prove that
$i^*L\in\ogenu {i^*H}{}{-M,M}$. 

Next we want to analyze the object $i^*L$, and for
this consider
the commutative square of open immersions
\[\xymatrix{
U\cap V \ar[rr]^-g\ar[d]_f & & V\ar[d]^k \\
U \ar[rr]_-i & & X
}\]
This square is cartesian, and flat base-change yields that
the following square of derived functors commutes
up to natural isomorphism
\[\xymatrix@C+30pt{
\Dqc(U\cap V)\ar[d]_{f_*} & \Dqc(V)\ar[d]^{k_*}\ar[l]_-{g^*} \\
\Dqc(U)  & \Dqc(X)\ar[l]^-{i^*}
}\]
Now: we are given in the category $\Dqcs {Z}(X)$ the distinguished
triangle $L\la F\la k_*k^*F\la\T L$, which the functor
$i^*$ takes to the distinguished triangle 
$i^*L\la i^*F\la i^*k_*k^*F\la\T i^*L$.
The isomorphism $i^*k_*\cong f_*g^*$ allows us to express
$i^*k_*k^*F$ as $f_*g^*k^*F$, and we have shown that our
triangle $i^*L\la i^*F\la i^*k_*k^*F\la\T i^*L$ is
isomorphic to some
triangle $i^*L\la i^*F\la f_*g^*k^*F\la\T i^*L$.
But now $i^*L\in\Dqcs{Z_1}(U)$, and 
\begin{eqnarray*}
  \Hom\big(\Dqcs{Z_1}(U)\,\,,\,\,f_*g^*k^*F\big)&\cong&\Hom\big(f^*\Dqcs{Z_1}(U)\,\,,\,\,g^*k^*F\big) \\
      &=& 0
\end{eqnarray*}
where the first isomorphism is by adjunction, while the vanishing
comes from $f^*\Dqcs{Z_1}(U)=0$. Hence the triangle must be
the one that comes from Bousfield localization: with $\ph:(U-Z_1)\la U$
the open immersion, the triangle
above is canonically
isomorphic to the
triangle obtained when we complete the unit
$\eta:i^*F\la \ph_*\ph^*i^*F$ of the adjunction
$\xymatrix{\ph^*\colon\Dqc(U)\ar@<0.5ex>[r] &
\Dqc(U-Z_1)\ar@<0.5ex>[l]\colon \ph_*
}$
to a triangle $L'\la i^*F\la \ph_*\ph^*i^*F\la\T L'$.

Thus we are down to the affine situation: we are given an affine
scheme $U$, a closed subset
$Z_1$ with quasicompact complement, and a perfect complex
$F'=i^*F\in\Dqc(U)$. And, completing the unit of adjunction
$\eta:F'\la \ph_*\ph^*F'$ to a distinguished triangle
$L'\la F'\la \ph_*\ph^*F'\la\T L'$, it suffices to show the inclusion
$L'\in\ogenu {H}{}{-M,M}$, with $H$ any compact generator
for $\Dqcs{Z_1}(U)$.

Now Lemma~\ref{L28.999} proves the special case where $F'=\co_U^{}$.
That is: in the triangle $L\la\co_U^{}\la\ph_*\ph^*\co_U^{}\la\T L$
we have that $L\in\ogenu {H}{}{-M,M}$. Let $F'\in\Dqc(U)$
be arbitrary and consider the commutative diagram with triangles as rows
\[\xymatrix{
F'\otimes L\ar@{.>}[d]_\alpha\ar[r] & F'\otimes\co_U^{}\ar[d]^\beta\ar[r] &
   F'\otimes \ph_*\ph^*\co_U^{}\ar[r]\ar[d]^\gamma &
   F'\otimes\T L\ar@{.>}[d]^{\T\alpha}\\
L'\ar[r] &F'\ar[r] &\ph_*\ph^* F'\ar[r] &\T L'
}\]
This means: the canonical isomorphism $\beta:F'\otimes\co_U^{}\la F'$
is because $\co_U^{}$ is the unit of the tensor product $\otimes$, the canonical
isomorphism $\gamma:F'\otimes \ph_*\ph^*\co_U^{}\la \ph_*\ph^* F'$
is the projection formula, and it is
classical that the solid square commutes. The dotted
map $\alpha$ is
defined by extending to a morphism of triangles\footnote{A more careful
analysis, also classical, shows that the isomorphism
$\alpha$ is in fact canonical.
We omit this here---the interested reader can find exhaustive accounts in
the literature.}. Since
$\beta$ and $\gamma$ are isomorphisms so is $\alpha$, giving that,
for any $F'\in\Dqc(U)$, we have $L'\cong F'\otimes L$.
Now if $F'$ is a perfect complex on $U=\spec R$, then there
exists an integer $N>0$ with  $F'\in\genu{\co_U^{}}{}{-N,N}$, and
hence
\[
L'\cong L\otimes F'\in L\otimes\genu{\co_U^{}}{}{-N,N}\subset\genu{L}{}{-N,N}
\subset\ogenu{H}{}{-M-N,M+N}\ ,
\]
and we have completed the proof of (ii), in the case where
$X=U\cup V$ with $U$ affine.

Before going
on: whenever we may assume (ii)
has already been proved,  for the inclusions $Z\subset X\supset V$,
we can apply it to the case where $F=G\in\Dqcs Z(X)$ is a compact generator. We
learn that in the triangle $L\la G\la k_*k^*G\la\T L$ we
have that $L\in\ogenu H{}{-M,M}$ for some $M>0$. Of course:
$G$ is a compact generator for $\Dqcs Z(X)$, and
$H\in\Dqcs{Z_1}(X)\subset\Dqcs Z(X)$ is a perfect complex. Hence
there must exist an integer $N>0$ with $H\in\genu G{}{-N,N}$. But 
then $L$ belongs to $\ogenu H{}{-M,M}\subset\ogenu G{}{-M-N,M+N}$,
and the triangle $G\la k_*k^*G\la\T L$ exhibits $k_*k^*G$ as belonging
to
\[\ogenu G{}{-M-N-1,M+N+1}*\ogenu G{}{-M-N-1,M+N+1}=\ogenu G{}{-M-N-1,M+N+1}\ .
\]

Now for the proof of (i), again in the
special case where $X=U\cup V$ with $U$ affine.
Let $F\in\Dqcs{Z\cap V}(V)$ be any perfect complex,
and let $G\in\Dqcs Z(X)$ be a compact generator. By the paragraph
above we may choose an integer $A>0$ with $k_*k^*G\in\ogenu G{}{-A,A}$.
Corollary~\ref{C27.3}, which we may apply since we
already know parts (i) and (ii) of
Theorem~\ref{T27.1}, tells us that there exists a
perfect complex 
$\wt F\in\Dqcs Z(X)$ with $k^*\wt F\cong F\oplus\T F$. Since
$\wt F$ is a compact object in $\Dqcs Z(X)$ and $G$ is a compact
generator we may, after increasing the integer $A>0$
if necessary,
assume further that
$\wt F\in\genu G{}{-A,A}$. But then
$k_*F\oplus k_*\T F\cong k_*k^*\wt F$ belongs
to
\[
k_*k^*\genu G{}{-A,A}\sub\genu{k_*k^*G}{}{-A,A}\sub
\ogenu G{}{-2A,2A}\ .
\]
And $k_*F$, being a direct summand of an object in $\ogenu G{}{-2A,2A}$,
must belong to $\ogenu G{}{-2A,2A}$.
This completes the proof of (i) in the case
where $X=U\cup V$ with $U$ affine, concluding the
proof of the Lemma.
\eprf

Next we prove the affine case of Theorem~\ref{T27.1}(iv).

\lem{L28.4}
Let $X$ be an affine scheme, and let $Z\subset X$ be a closed subset
with quasicompact complement. Then the category $\Dqcs Z(X)$ is
weakly approximable.

More precisely: with the compact generator $B\in\Dqcs Z(X)$ as
in Construction~\ref{C28.2}, and given any object $F$ in the category
$\Dqcs Z(X)^{\leq0}$ for the standard \tstr, there exists
a triangle $E\la F\la D\la\T E$ in $\Dqcs Z(X)$ with
$E\in\ogenu B{}{0,0}$ and with $D\in\Dqcs Z(X)^{\leq-1}$.
\elem

\prf
Writing $X=\spec R$ and using the equivalence $\Dqc(X)\cong\D(R)$, we
have that the given object $F\in\Dqcs Z(X)\cong\D_Z^{}(R)$
is assumed to belong to
$\D_Z^{}(R)^{\leq 0}$. We may choose a free
$R$--module $L$ and a
surjection $L\la H^0(F)$, which can
be lifted to a morphism $\psi:L\la F$ in
$\D(R)$. And now every basis element $i$ of $L$ corresponds to
a morphism $i:R\la L$, and the composite
$R\stackrel i\la L\stackrel\psi\la F$ is a morphism from
$R$ to a complex
with cohomology supported in $Z$.
And, proceeding as in  
the proof of \cite[Corollary~5.7.1(i)]{Lipman-Neeman07}, one
shows that
the map $\psi\circ i:R\la F$ must factor as $R\la B_{n_i}\la F$
for some $n_i>0$ depending on $i$, and with
$B_{n_i}$ as in Construction~\ref{C28.2}.
Assembling all these maps, the morphism $\psi:L\la F$ can be
factored as
\[\xymatrix{
L\ar[r] &\ds\coprod_{i\in I} B_{n_i}\ar[r]^-\alpha &
F  
}\]
where $I$ is the set of basis elements of $L$. Now
setting $E=\coprod_{i\in I} B_{n_i}$
we have that $E$ belongs to $\ogenu B{}{0,0}$,
and the morphism $H^0(\alpha):H^0(E)\la H^0(F)$ is an
epimorphism. The long
exact sequence in cohomology teaches us that,
in the triangle $E\stackrel\alpha\la F\la D\la\T E$,
we must have $D\in\D_Z^{}(X)^{\leq-1}$.
\eprf

And the time has finally come for

\medskip

\nin
\emph{Proof of Theorem~\ref{T27.1}(iv).}\ \ Lemma~\ref{L28.4} proves
the theorem for affine $X$ and, as in the proof of
parts (i), (ii) and (iii) of Theorem~\ref{T27.1}
back at the end of Section~\ref{S27}, the general case
follows by induction on the number of open affines needed to
cover $X$. We may again assume that
$X=U\cup V$, with $U$ affine and with the theorem true for $V$.
In the notation developed back then,
we constructed a perfect complex $G\in\Dqcs Z(X)$ and a compact
generator $S\in\Dqcs{Z_1}(X)$, in such a way that $S$ is a direct
summand of $G$ and $k^*G\in\Dqcs{Z\cap V}(V)$ is a compact generator.
It then formally followed that $G$ is a compact generator for
$\Dqcs Z(X)$. Let $S$ and $G$ be as back then.

Lemma~\ref{L28.1000}(i) now tells us
  that we may choose an integer $M>0$ with
$k_*k^*G\in\ogenu G{}{-M,M}$.
Next choose an
integer $N>0$ such that
$k_*\Dqcs {Z\cap V}(V)^{\leq0}\subset\Dqcs Z(X)^{\leq N}$. Increasing
the integer $N$ if necessary, we can furthermore assume
$\Hom\left(G,\Dqcs Z(X)^{\leq-N}\right)=0$.
For the affine scheme $U$, Lemma~\ref{L28.4} tells
us that $\Dqcs{Z_1}(U)\cong\Dqcs{Z_1}(X)$ is weakly approximable. And
we are assuming that $\Dqcs{Z\cap V}(V)$ is already known
to be weakly approximable. Increasing the integer $N>0$ some more
gives
\begin{itemize}
\item
For any objects
$F_1\in\Dqcs{Z_1}(X)^{\leq0}$ and $F_2\in\Dqcs{Z\cap V}(V)^{\leq0}$ there
exist triangles 
\[\xymatrix@R-20pt{
E_1\ar[r] & F_1\ar[r] & D_1\ar[r] & \T E_1 \\
E_2\ar[r] & F_2\ar[r] & D_2\ar[r] & \T E_2
}\]
with $E_1\in\ogenu{S}{}{-N,N}\subset\ogenu{G}{}{-N,N}$,
with $E_2\in\ogenu{k^*G}{}{-N,N}$,
with $D_1\in\Dqcs{Z_1}(X)^{\leq-1}$ and with
$D_2\in\Dqcs{Z_V}(V)^{\leq-1}$.
\end{itemize}

We have now chosen all the integers given by the inductive hypothesis
and by the Lemmas so far. Now take any object $F\in\Dqcs Z(X)^{\leq0}$.
Then $k^*F$ belongs to $\Dqcs{Z_V}(V)^{\leq0}$, and by combining
$\bullet$ above with \cite[2.2.1]{Neeman17A} we may choose
in $\Dqcs{Z_V}(V)$ a distinguished triangle
$E'\la k^*F\la D'\la\T E'$, with $E'\in\ogenu{k^*G}{}{-2N,N}$ and with
$D'\in\Dqcs{Z_V}(V)^{\leq-N-1}$. And now the triangle
$k_*E'\la k_*k^*F\la k_*D'\la\T k_*E'$ satisfies
\[
k_*D'\quad\in\quad k_*\Dqcs{Z_V}(V)^{\leq-N-1}\sub\Dqcs Z(X)^{\leq-1}
\]
and
\[
k_*E'\quad\in\quad k_*\ogenu{k^*G}{}{-2N,N}\sub\ogenu{k_*k^*G}{}{-2N,N}
\sub\ogenu{G}{}{-M-2N,M+N}
\]
where the last inclusion comes from the fact that $k_*k^*G$
belongs to $\ogenu G{}{-M,M}$.

Now complete in $\Dqcs Z(X)$ the unit of adjunction $\eta:F\la k_*k^*F$
to a distinguished triangle $L\la F\la k_*k^*F\la \T L$.
Since $F$ belongs to $\Dqcs Z(X)^{\leq0}$ we deduce first that
$k^*F\in\Dqcs{Z_V}(V)^{\leq0}$, and then
that $k_*k^*F\in\Dqcs Z(X)^{\leq N}$.
Now the triangle $k_*k^*F\la\T L\la\T F$ exhibits
$\T L$ as the extension of two objects in $\Dqcs Z(X)^{\leq N}$,
hence $\T L\in\Dqcs Z(X)^{\leq N}$. But $\T L$ is annihilated by
$k^*$, and must therefore belong to
$\Dqcs{Z_1}(X)^{\leq N}$. We may therefore choose in
$\Dqcs{Z_1}(X)$ an exact triangle $E''\la \T L\la D''\la\T E''$
with $E''\in\ogenu{G}{}{-M-4N,2N}$
and with
$D''\in\Dqcs{Z_1}(X)^{\leq-M-3N-1}\subset\Dqcs Z(X)^{\leq-M-3N-1}$.

Now consider the diagram where the columns are triangles
\[\xymatrix{
  k_*E'\ar[d] & E''\ar[d]\\
  k_*k^* F' \ar[r]\ar[d] & \T L\ar[d]\\
  k_* D' & D''  
}\]
The composite from top left to bottom right is a morphism
from the object $k_*E'\in\ogenu G{}{-M-2N,M+N}$ 
to the object $D''\in\Dqcs Z(X)^{\leq -M-3N-1}$. But $N$ was chosen
so that $\Hom\left(G,\Dqcs Z(X)^{\leq-N}\right)=0$, hence
$\Hom\left(-,\Dqcs Z(X)^{\leq-M-3N-1}\right)$ vanishes on all
$\T^iG$ with $i\leq M+2N+1$, and therefore
$\Hom\left(-,\Dqcs Z(X)^{\leq-M-3N-1}\right)$ also vanishes on
all of $\ogenu G{}{-M-2N,M+N}$. In particular it vanishes on
$k_*E'$. Thus the composite $k_*E'\la D''$ vanishes, and we deduce
that there must exist a map $k_*E'\la E''$ rendering commutative
the square
\[\xymatrix{
k_*E'\ar[r]\ar[d] & E''\ar[d] \\
k_*k^*F\ar[r] &\T L
}\]
And now we may complete to a $3\times3$ diagram where the rows and
columns are triangles
\[\xymatrix{
\Tm E''\ar[r]\ar[d] & E\ar[r]\ar[d] &k_*E'\ar[r]\ar[d] & E''\ar[d] \\
L\ar[r]\ar[d] & F\ar[r]\ar[d] &k_*k^*F\ar[d]\ar[r] &\T L\ar[d] \\
\Tm D''\ar[r] &D\ar[r] &k_*D'\ar[r] & D''
}\]
The fact that $k_*D'\in\Dqcs Z(X)^{\leq-1}$ and
$\Tm D''\in\Dqcs Z(X)^{\leq-M-3N}$ certainly forces
$D\in\Dqcs Z(X)^{\leq-1}$. And as
$k_*E\in\ogenu G{}{-M-2N,M+N}$ while
$\Tm E''\in\ogenu{G}{}{-M-4N+1,2N+1}$, this combines
to give that
$E$ must lie in $\ogenu{G}{}{-M-4N,M+4N}$.
The triangle $E\la F\la D$ now proves the weak approximability
of $\Dqcs Z(X)$.\hfill{$\Box$}

\section{More about approximable and weakly approximable triangulated
categories}
\label{S400}

Let $X$ be a quasicompact, quasiseparated scheme, and let $Z\subset X$
be a closed subset with quasicompact complement.
Theorem~\ref{T27.1} proves that (1) the category $\ct=\Dqcs Z(X)$ has a single compact
generator, (2) the standard \tstr\ is in the preferred
equivalence class, and (3) the
category is weakly approximable. From (2) it follows that the preferred
subcategories $\ct^b$, $\ct^-$ and $\ct^+$,
introduced in \cite{Neeman17A}, agree (respectively) with 
the classical $\Dqcsb Z(X)$, $\Dqcsmi Z(X)$ and $\Dqcspl Z(X)$.
 
And, in combination with
Theorem~\ref{T27.1}, 
Theorem~\ref{T29.1} shows that, 
with respect to the preferred equivalence
class of \tstr s, the pseudocoherent complexes on $X$,
with cohomology supported on $Z$, are precisely
the objects
of $\ct=\Dqcs Z(X)$ which can be approximated arbitrarily well by
objects in $\ct^c=\dperfs ZX$. Thus the category $\ct^-_c$ of
\cite{Neeman17A} agrees with the category of pseudocoherent complexes with
cohomology supported on $Z$. And $\ct^b_c=\ct^-_c\cap\ct^b$ identifies with
the pseudocoherent complexes whose cohomology vanishes in all but finitely
many degrees. 

If $X$ is noetherian this simplifies to
$\ct^-_c=\dmcohs Z(X)$ and 
$\ct^b_c=\dcohs Z(X)$.

Until this article these results were only known when $Z=X$ and the scheme
is assumed separated. The proof used the fact that the category $\Dqc(X)$ 
is approximable, and a theorem
saying that the proof of approximability delivers
information about the preferred
equivalence class of \tstr s.

\rmk{R400.1}
In Theorem~\ref{T27.1}(iv) we 
proved the weak approximability of $\Dqcs Z(X)$. Now
we will show that $\Dqcs Z(X)$ is not approximable in general.
For simplicity let us assume that $X=\spec R$ is affine, and let
$B$ be the compact generator of $\Dqcs Z(X)$
given in Construction~\ref{C28.2}, that is $B$ is the
complex
\[
B\eq\bigotimes_{i=1}^r\left(R\stackrel{f_i}\la R\right)\ .
\]
Let $I\lhd R$ be the ideal generated by $\{f_1,f_2,\ldots,f_r\}$.

The object $B$ is annihilated by every element of the ideal $I$,
and the objects of $\Coprod_1^{}\big(B(-\infty,\infty)\big)$ are
coproducts of suspensions of $B$. Hence the ideal $I$ annihilates
every object of $\Coprod_1^{}\big(B(-\infty,\infty)\big)$.
Next: the inductive formula tells us that
\[
\Coprod_{n+1}^{}\big(B(-\infty,\infty)\big)\eq
\Coprod_1^{}\big(B(-\infty,\infty)\big)*\Coprod_n^{}\big(B(-\infty,\infty)\big)\ .
\]
Hence induction shows that the ideal $I^n$ annihilates
$\Coprod_n^{}\big(B(-\infty,\infty)\big)$. Consequently, for every
object $E\in\Coprod_n^{}\big(B(-\infty,\infty)\big)$ we must have
that $I^nH^0(E)=0$.
If $E\la F\la D$ is a triangle with $E\in\Coprod_n^{}\big(B(-\infty,\infty)\big)$,
with $F\in D(R)^{\leq0}$ and with $D\in\D(R)^{\leq-1}$, then
the exact sequence $H^0(E)\la H^0(F)\la H^0(D)=0$
would tell us that $I^nH^0(F)=0$. And this does not normally
happen if $F=R/I^{n+1}$; it happens exactly for those ideals with 
$I^n=I^{n+1}$. This does happen, but rarely: an example is $I=0$,
in which case $Z=X$.

Thus we cannot usually find an integer $n$ such that, for every 
$F\in\Dqcs Z(X)^{\leq0}$, there exists a triangle $E\la F\la D$ with
$E\in\Coprod_n^{}\big(B[-n,n]\big)$ and with $D\in\Dqcs Z(X)^{\leq-1}$.
\ermk 

\rmk{R400.3}
Theorem~\ref{T27.1}(iv) shows that the category $\Dqcs Z(X)$ is
always weakly approximable, and in
Remark~\ref{R400.1} we showed that we can only
expect it to be approximable for very special $Z\subset X$.

Of course:
if $Z=X$ and $X$ is separated then we know, from
the examples in \cite[Section~3]{Neeman17A}, that
the category
$\Dqcs Z(X)$ is in fact approximable. But this isn't an easy theorem,
or more accurately the current proof isn't easy.

It is relevant to note that the very general weak
approximability result, of the current article,  was proved
by reduction to the case where $X$ is affine. By contrast the known
proof, that $\Dqc(X)$ is approximable for separated $X$, is by reduction
to the case where $X$ is projective. Reducing to the projective
case is technically more involved, and forces on us the assumption
that $X$ is separated.

Let's come back to 
the current article: in a precursor version 
we made restrictive assumptions that
permitted us to apply the theory of approximable triangulated categories
as it stood,
unmodified and unextended. Section~\ref{S27} evolved 
in order to prove the theorems of this article,
about bounded \tstr s on
$\dperfs ZX$ and on $\dcohs Z(X)$, in the current, much extended
generality.
And Section~\ref{S28} was prompted by an application
unrelated to anything in the current article: the weak approximability
of $\Dqcs Z(X)$ turns out to be relevant to a question
studied by Canonaco and Stellari.

The computation of the intrinsic $\ct^-$, $\ct^+$, $\ct^b$, 
$\ct^-_c$ and $\ct^b_c$ of \cite{Neeman17A}, for $\ct=\Dqcs Z(X)$,
came as an unexpected byproduct.

The theory of approximable triangulated categories should be viewed as
work-in-progress, with extensions and modifications encouraged. The
theory should evolve to apply more widely, with this article being
a manifestation. 
\ermk

Needless to say: it remains to explore the wider implications of
the results developed here, to other approximable or nearly-approximable
triangulated categories.

\providecommand{\bysame}{\leavevmode\hbox to3em{\hrulefill}\thinspace}
\providecommand{\MR}{\relax\ifhmode\unskip\space\fi MR }
\providecommand{\MRhref}[2]{%
  \href{http://www.ams.org/mathscinet-getitem?mr=#1}{#2}
}
\providecommand{\href}[2]{#2}

\end{document}